\documentclass{article}[letterpaper, 14pt]
\usepackage{graphicx}
\usepackage{amsmath,bm}
\usepackage{amsmath}
\usepackage{graphicx}
\usepackage{amsfonts}
\usepackage{xcolor}
\usepackage[letterpaper, total={6in, 8in}]{geometry}
\usepackage{pdfpages}

\usepackage{authblk}

\title{Representation of tensor functions using low-order structural tensor set: two-dimensional point groups}
\author[1]{Mohammad Madadi\thanks{Corresponding author e-mail: mmadadi1@binghamton.edu}}
\author[2]{Lin Cheng}
\author[1]{Pu Zhang\thanks{Corresponding author e-mail: pzhang@binghamton.edu}}
\affil[1]{\textit{Department of Mechanical Engineering, State University of New York at Binghamton, Binghamton, NY 13902, USA}}
\affil[2]{\textit{Department of Mechanical Engineering, University of Maryland, College Park, MD 20742, USA}}

\begin{document}

\maketitle

\begin{abstract}
The representation theory of tensor functions is essential to constitutive modeling of materials including both mechanical and physical behaviors. Generally, material symmetry is incorporated in the tensor functions through a structural or anisotropic tensor that characterizes the corresponding point group. The general mathematical framework was well-established in the 1990s. Nevertheless, the traditional theory suffers from a grand challenge that many point groups involve fourth or sixth order structural tensors that hinder its practical applications in engineering. 
Recently, researchers have reformulated the representation theory and opened up opportunities to model anisotropic materials using low-order (i.e., 2nd- order and lower) structural tensors only, although the theory was not fully established. This work aims to fully establish the reformulated representation theory of tensor functions for all two-dimensional point groups. It was found that each point group needs a structural tensor set to characterize the symmetry. For each two-dimensional point group, the structural tensor set is proposed and the general tensor functions are derived. Only low-order structural tensors are introduced so researchers can readily apply these tensor functions for their modeling applications. The theory presented here is useful for constitutive modeling of materials in general, especially for composites, nanomaterials, soft tissues, etc. 
\end{abstract}

\section{Introduction}
Constitutive modeling of materials \cite{ottosen2005mechanics} plays a crucial role in continuum mechanics and materials science. Contemporary constitutive modeling encompasses not only mechanical behaviors (e.g., stress-strain relations,
yield surface, and failure criterion) but also physical properties (e.g., conductivity, dielectric, etc.) and mechano-physical behaviors (i.e., strain-dependent physical properties). These mechanical behaviors or physical properties are usually described by scalar- or tensor-valued functions with multiple arguments including field and state variables. Generally, most materials are anisotropic and they fall into one of the 12 two-dimensional (2D) and 39 three-dimensional (3D) point groups \cite{zheng_theory_1994}, if both crystallographic and continuous groups are counted. The representation theory of tensor functions offers general forms of constitutive models that satisfy both frame-indifference and material symmetry principles. It is among the most significant theories in continuum mechanics and lays the foundation for constitutive modeling. The establishment of this theory dated back to the mid-20th century. Researchers like Rivlin (1948) \cite{rivlin_large_1948}, Rivlin and Ericksen (1955) \cite{rivlin_stress_1997}, Rivlin (1955) \cite{rivlin_further_1955}, Pipkin and Rivlin (1959) \cite{pipkin_formulation_1959}, and Noll (1970) \cite{noll_representations_1970} introduced complete and irreducible polynomial representations of isotropic tensor functions. Moreover, Wang \cite{wang_representations_1969, wangb_representations_1969, wang_new_1970, wang_corrigendum_1971}, Smith \cite{smith_fundamental_1970, smith_isotropic_1971}, and Boehler \cite{boehler_irreducible_1977} derived the first general relations for isotropic scalar-, vector-, 2nd-order symmetric, and antisymmetric tensor-valued functions of any finite number of 2nd-order symmetric and skew-symmetric tensors, and vectors in 3D space. Afterwards, the late 20th century witnessed further development on the representation theory of anisotropic tensor functions. Specifically, Boehler, Spencer, Liu and Betten \cite{boehler_lois_1978, boehler_simple_1979, Boehler1987, spencer_formulation_1982, liu_representations_1982, betten_irreducible_1987} developed a method of isotropic extension by introducing structural tensors (a.k.a., anisotropic tensors). The structural tensors characterize material symmetry and are invariant under any symmetry operation in the corresponding point group. This approach converts an anisotropic tensor function into an extended isotropic tensor function by including additional structural tensor arguments. The structural tensors of all 2D and 3D point groups can be found in works of Zheng and Xiao \cite{zheng_two-dimensional_1997, zheng_note_1994, zheng_theory_1994, xiao_isotropic_1996}. Researchers can refer to a comprehensive review paper by Zheng \cite{zheng_theory_1994} for more details. 

Although this representation theory of anisotropic tensor functions is elegant, its practical application is somehow limited and the development is stalled. On the one hand, it does not provide specific functions for constitutive models, but rather the general mathematical formulation. One has to figure out the specific functions either empirically or from a trial-and-error process. On the other hand, a grand challenge is that most point groups exhibit higher-order structural tensors (i.e., $3^{rd}$ order and higher) that impede practical modeling usage \cite{zheng_theory_1994}. For example, a $6^{th}$ order structural tensor was defined for the group $\mathcal{C}_{6v}$ and it is very hard to devise specific tensor functions including this term. Generally, constitutive modeling with structural tensors of order higher than two is complicated and impractical. Among all the material symmetries, besides orthotropic and transversely isotropic, only point groups in triclinic, monoclinic, and rhombic crystal systems can be characterized with lower-order structural tensors (i.e., $2^{nd}$ order and lower). For other point groups, very little constitutive modeling work has been done yet.  

In order to circumvent higher-order structural tensors, Man and Goddard \cite{man2018remarks} reformulated the representation theory of tensor functions in 2018, which sheds some light to develop constitutive models for practical usage. In the original formulation developed by Boehler and Liu \cite{boehler_simple_1979, liu_representations_1982}, the structural tensors must be invariant under any symmetry operation in the point group. In contrast, the reformulation by Man and Goddard \cite{man2018remarks} relaxed this strict requirement for structural tensors but added additional symmetry constraints on the tensor functions afterwards.  Man and Goddard's work opened a new era in the representation of anisotropic tensor functions, allowing constitutive modeling for all point groups with lower-order structural tensors only. In their paper, Man and Goddard provided a couple of examples to introduce how to use the reformulation. There is still a large knowledge gap on the reformulated representation theory and its applications in all 2D and 3D point groups. This work aims to address this need, establish the reformulated representation theory for all 2D point groups, and partially fill this knowledge gap.  

In the following sections, we first revisit the representation theory of tensor functions in Section \ref{sec:rep} to prepare readers for the background theory. In Section \ref{sec:struct_tensor_set}, we introduce a concept "structural tensor set" that is the key to this reformulated representation and present the structural tensor sets for all 2D point groups. Note that six 2D point groups ($\mathcal{C}_1$, $\mathcal{C}_{1v}$, $\mathcal{C}_2$, $\mathcal{C}_{2v}$, $\mathcal{C}_\infty$, $\mathcal{C}_{\infty v}$) have lower-order structural tensors so their representation theory can adopt the original formulation of Boehler and Liu, as listed in Section \ref{sec:rep}. In contrast, for the six 2D point groups ($\mathcal{C}_4$, $\mathcal{C}_{4v}$, $\mathcal{C}_3$, $\mathcal{C}_{3v}$, $\mathcal{C}_6$, $\mathcal{C}_{6v}$) with higher-order structural tensors, the reformulation of Man and Goddard has to be employed to establish the representation theory. Hence, we present the reformulated representation theory for these six 2D point groups in Sections \ref{sec:C4}-\ref{sec:C6}. It is worth mentioning that this paper uses the Schoenflies notation. For more details on the various notation systems, readers are encouraged to check the SI and consult \cite{de2012structure}. Finally, we add more discussions and remarks in Section \ref{sec:remarks} regarding the reformulated representation theory and its applications. Additionally, it is worth noting that we have listed the symmetry operations and multiplication tables for all 2D point groups in the Supplementary Information. To our knowledge, these multiplication tables have not been published elsewhere yet. The Supplementary Information is not only useful for the theories outlined in this work but also to the broader research community. 

\section{Introduction to representation theory of tensor functions} \label{sec:rep}
In this section, a brief introduction to the representation theory of tensor functions is provided, including both isotropic and anisotropic tensor functions. In addition, tensor generators and functional bases of point groups with lower-order structural tensors ($\mathcal{C}_1$, $\mathcal{C}_{1v}$, $\mathcal{C}_2$, $\mathcal{C}_{2v}$, $\mathcal{C}_\infty$, $\mathcal{C}_{\infty v}$) are presented for users to generate tensor functions. At the end, the reformulated representation theory of Man and Goddard is introduced. Although we focus on tensor functions in the 2D space, the theories and equations can be generalized to 3D readily. 

\subsection{Isotropic tensor functions}

A systematic and in-depth introduction to the representation theory of isotropic tensor functions can be found in \cite{Itskov2007, zheng_note_1994, zheng_two-dimensional_1997}. Below is a very brief summary to the key concepts and equations.

Let $\mathcal{O}(2)$ denote the 2D orthogonal group, $\mathbf{Q}$ be an orthogonal transformation in 2D. A scalar-valued function $\psi(\mathbf v, \mathbf{A}, \mathbf{W})$ and a 2nd-order symmetric tensor-valued function $\mathbf{T}(\mathbf v, \mathbf{A}, \mathbf{W})$ are called isotropic if
\begin{align} \label{eq 2.1}
	\begin{array}{ll}
		\psi(\mathbf v, \mathbf{A}, \mathbf{W})=\psi(\mathbf{Q}\mathbf{v}, \mathbf{Q}\mathbf{A}			\mathbf{Q}^T, \mathbf{Q}\mathbf{W}\mathbf{Q}^T)\\
		\mathbf{Q}\mathbf{T}(\mathbf{v}, \mathbf{A}, \mathbf{W})\mathbf{Q}^T=\mathbf{T}(\mathbf{Q}		\mathbf{v}, \mathbf{Q}\mathbf{A}\mathbf{Q}^T,  \mathbf{Q}\mathbf{W} \mathbf{Q}^T) 
	\end{array}   ;\forall \mathbf{Q} \in \mathcal{O}(2)
\end{align}
where $\mathbf v$, $\textbf{A}$, and $\mathbf{W}$ are sets of vectors, 2nd-order symmetric and skew-symmetric tensors, respectively. In continuum mechanics, Eq. (\ref{eq 2.1}) guarantees the frame-indifference. 

The representation theory for an isotropic scalar-valued function states that $\psi(\mathbf v, \mathbf{A}, \mathbf{W})$ can be expressed as a function of invariants $I_k$ of its arguments $\mathbf v$, $\textbf{A}$, and $\mathbf{W}$,
\begin{align} \label{eq 2.2}
	\psi(\mathbf v, \mathbf{A}, \mathbf{W})=\psi(I_k)
\end{align}
where the complete set of invariants $I_k(k=1,2,...,r)$ is called functional basis (or integrity basis). A basis is considered minimal if it contains only irreducible invariants \cite{korsgaard_representation_1990}. Once the arguments $\mathbf v$, $\textbf{A}$, and $\mathbf{W}$ are provided, one can follow equations in \cite{korsgaard_representation_1990, Boehler1987} to find the invariants $I_k$. 

The representation of an isotropic 2nd-order symmetric tensor-valued function $\mathbf{T}(\mathbf{v}, \mathbf{A}, \mathbf{W})$ is expressed as a linear combination of a generating set of tensors $\mathbf G_i$. In general, $\mathbf{T}(\mathbf{v}, \mathbf{A}, \mathbf{W})$ follows the following form.  
\begin{align} \label{eq 2.3}
	\begin{array} {ll}
		\mathbf{T}(\mathbf{v}, \mathbf{A}, \mathbf{W})=\sum\alpha_i \mathbf G_i \\
		\alpha_i = \alpha_i(I_1,I_2,...,I_r)
	\end{array}
\end{align}
where $\alpha_i$ is a scalar coefficient function of the invariants $I_k$, $\mathbf G_i$ is a tensor generator in the generating set. Once the arguments $\mathbf v$, $\textbf{A}$, and $\mathbf{W}$ are provided, one can find the tensor generating set $\mathbf{G}_i$ following methods and formulae presented in \cite{Boehler1987, zheng_theory_1994, korsgaard_representation_1990}.

Most studies on the representation of 2D isotropic functions have focused on their relation to 3D functions \cite{adkins_symmetry_1959, adkins_further_1960, ajm_theory_1971}. They are deduced from the representation of 3D isotropic functions by reducing the Cayley-Hamilton theorem for two-dimensional tensors. However, this procedure is not only tedious and laborious but also limited to the polynomial representation (not the general form). In 1990, Korsgaard \cite{korsgaard_representation_1990} provided functional basis and tensor generators for the scalar, vector, second-order symmetric, and skew-symmetric tensor-valued isotropic functions in 2D through the direct method. Thus, for a scalar- or tensor-valued isotropic function $\psi(\mathbf v, \mathbf{A},\mathbf{W})$ and $\mathbf T(\mathbf v, \mathbf{A},\mathbf{W})$, the functional basis and generating set can be obtained from Table \ref{Table 1} and Table \ref{Table 2}, respectively.

For applications in constitutive modeling of isotropic materials, the isotropic scalar-valued function $\psi(\mathbf v, \mathbf{A},\mathbf{W})$ can be used to describe strain energy, free energy, yield surface, failure criterion, etc. On the other hand, the isotropic $2^{nd}$-order symmetric tensor-valued function $\mathbf T(\mathbf v, \mathbf{A},\mathbf{W})$ can be used to describe stress, conductivity, dielectric property, etc. In addition, it will be shown below that the representation of anisotropic functions is also built upon the isotropic theory.

\begin{table}[!h]
	\caption{Functional basis of  $\psi(\mathbf v, \mathbf{A},\mathbf{W})$ in 2D space \cite{korsgaard_representation_1990}}
	\label{Table 1}
	\begin{tabular}{ll}
		\hline
		Functional basis &  \\
		\hline
		$tr\mathbf{A}_i$, $tr\mathbf{A}_i^2$ &  \\
		$tr\mathbf{A}_i\mathbf{A}_j$& $i,j= 1,2, ..., N; i<j$ \\
		$\mathbf{v}_m . \mathbf{v}_m$ &	$m,n= 1,2, ..., M; m<n$\\
		$\mathbf{v}_m . \mathbf{v}_n$ & $p,q= 1,2, ..., P; p<q$	\\
		$tr\mathbf{W}_p^2$ &		\\
		$tr\mathbf{W}_p\mathbf{W}_q$ &	\\
		$\mathbf{v}_m .\mathbf{A}_i \mathbf{v}_m$ &		\\
		$\mathbf{v}_m .\mathbf{A}_i \mathbf{v}_n$ &		\\
		$\mathbf{v}_m .\mathbf{W}_p \mathbf{v}_n$ &		\\
		$tr\mathbf{A}_i\mathbf{A}_j\mathbf{W}_p$ &		\\
		$\mathbf{A}_i\mathbf{v}_m .\mathbf{W}_p \mathbf{v}_m$ &		\\\hline
	\end{tabular}
	\vspace*{-4pt}
\end{table} 

\begin{table}[!h]
	\caption{Tensor generators of $\mathbf{T}(\mathbf v, \mathbf{A},\mathbf{W})$ in 2D space \cite{korsgaard_representation_1990}}
	\label{Table 2}
	\begin{tabular}{ll}
		\hline
		Tensor generators &  \\
		\hline
		$\mathbf{I}$\footnotemark[1] &  \\
		$\mathbf{A}_i$ & $i,j= 1,2, ..., N$ \\
		$\mathbf{v}_m \otimes \mathbf{v}_m$ &	$m,n= 1,2, ..., M; m<n$\\
		$\mathbf{v}_m \otimes \mathbf{v}_n +\mathbf{v}_n \otimes \mathbf{v}_m $ & $p,q= 1,2, ..., P$	\\
		$\mathbf{v}_m \otimes \mathbf{W}_p \mathbf{v}_m +\mathbf{W}_p  \mathbf{v}_m \otimes \mathbf{v}_m $ \quad &		\\
		$\mathbf{A}_i \mathbf{W}_p - \mathbf{W}_p \mathbf{A}_i$ &		\\\hline
	\end{tabular}
	\footnotetext[1]{$\mathbf I$ is the second-order identity tensor in 2D space.}
	\vspace*{-4pt}
\end{table} 

\subsection{Anisotropic tensor functions: original formulation by Boehler and Liu}
Anisotropic tensor functions are needed for constitutive modeling of anisotropic materials. The traditional theory was developed by Boehler \cite{Boehler1987} and Liu \cite{liu_representations_1982} and reviewed by Zheng \cite{zheng_theory_1994}. The theory is briefly summarized below. Note that in this study, we limit the theory to 2D space.

The scalar- and tensor-valued functions $\psi(\mathbf v, \mathbf{A},\mathbf{W})$ and $\mathbf T (\mathbf v, \mathbf{A},\mathbf{W})$ are anisotropic if the condition in (\ref{eq 2.1}) is not satisfied $\forall \mathbf{Q} \in \mathcal{O}(2)$. In other words, considering $\mathcal{G}$ the point group of interest, the scalar- and tensor-valued functions are called anisotropic if 
\begin{align} \label{eq 2.4}
	\begin{array}{ll}
		\psi(\mathbf v, \mathbf{A}, \mathbf{W})=\psi(\mathbf{Q}\mathbf{v}, \mathbf{Q}\mathbf{A}			\mathbf{Q}^T, \mathbf{Q}\mathbf{W}\mathbf{Q}^T)\\
		\mathbf{Q}\mathbf{T}(\mathbf{v}, \mathbf{A}, \mathbf{W})\mathbf{Q}^T=\mathbf{T}(\mathbf{Q}		\mathbf{v}, \mathbf{Q}\mathbf{A}\mathbf{Q}^T, \mathbf{Q}\mathbf{W} \mathbf{Q}^T) 
	\end{array}    ;\forall \mathbf{Q} \in \mathcal{G}
\end{align} 
Boehler and Liu proposed a method of isotropic extension to derive the representation theory of anisotropic tensor functions. The idea is to introduce structural tensors $\mathbf{M}$ in the arguments and extend the anisotropic tensor functions to isotropic ones, and then use Eqs. (\ref{eq 2.2}) and (\ref{eq 2.3}) to derive the representation. Before presenting the theory, we will introduce the structural tensors below first.

One can either define a single or multiple structural tensors for each point group \cite{zheng_theory_1994}. Each structural tensor is invariant under any symmetry operation $\mathbf{Q}$ in the point group $\mathcal{G}$. The structural tensors $\mathbf{M}$ can include vectors, second-, and higher-order tensors that characterize the material symmetry. 
The definition of structural tensors initially presented in \cite{ lokhin_nonlinear_1963, smith_anisotropic_1957,boehler_simple_1979,boehler_lois_1978, zheng_theory_1994} is introduced as follows. 
For a point group $\mathcal{G}$, if there exists a set of tensors $\mathbf M =\{\mathbf M_1, \mathbf M_2,...,\mathbf M_s\}$ such that any orthogonal transformation $\mathbf Q \in \mathcal O (2)$ belongs to $\mathcal G$ if and only if each tensor $\mathbf{M}_i$ is invariant, we say that tensors $\mathbf M_1$, $\mathbf M_2$,...,$\mathbf M_s$ characterize $\mathcal G$ and call them structural tensors. Mathematically, the definition is 
\begin{align} \label{eq 2.6}
	\mathcal G = \{ \mathbf Q \in \mathcal O (2)| \langle \mathbf Q\rangle\mathbf M_i = \mathbf M_i; i=1,2,...,s\}
\end{align}
where the orthogonal transformation operator $\langle \mathbf Q\rangle$ is defined after Zheng \cite{zheng_theory_1994}, as:
\begin{align} \label{eq 2.7}
	\begin{array} {ll}
		\langle\mathbf Q \rangle \mathbf v = \mathbf Q \mathbf v=  \mathbf Q_{ij} \mathbf v_j \\
		\langle\mathbf Q \rangle \mathbf A = \mathbf{Q} \mathbf{A} \mathbf{Q}^T= \mathbf Q_{ip} \mathbf Q_{jq} \mathbf A_{pq}\\ 
		\langle\mathbf Q \rangle \mathbb A = \mathbf Q_{ip} \mathbf Q_{jq} \mathbf Q_{kr}...\mathbf Q_{mt} \mathbb A_{pqr...t}
	\end{array}
\end{align}
where $\mathbf v$, $\mathbf A $, and $\mathbb A$ are first-, second-, and higher-order tensors, respectively.

By adding structural tensors $\mathbf M$ in the arguments, we can obtain an isotropic extension of the anisotropic scalar- or tensor-valued functions  $\hat{\psi}(\mathbf v, \mathbf{A}, \mathbf{W}, \mathbf M)$  and  $\mathbf{\hat T}(\mathbf v, \mathbf{A}, \mathbf{W}, \mathbf M)$ as follows
\begin{align} \label{2.5}
	\begin{array} {ll}
		\psi(\mathbf v, \mathbf{A}, \mathbf{W})=\hat{\psi}(\mathbf v, \mathbf{A}, \mathbf{W}, \mathbf M) \\
		\mathbf{T}(\mathbf v, \mathbf{A}, \mathbf{W})=\mathbf{\hat T}(\mathbf v, \mathbf{A}, \mathbf{W}, \mathbf M) 
	\end{array}
\end{align}
and then $\hat{\psi}(\mathbf v, \mathbf{A}, \mathbf{W}, \mathbf M)$  and  $\mathbf{\hat T}(\mathbf v, \mathbf{A}, \mathbf{W}, \mathbf M)$ are isotropic tensor functions satisfying the condition (\ref{eq 2.1}), i.e.
\begin{align} \label{2.5 prime}
	\begin{array} {ll}
		{\hat \psi}(\mathbf v, \mathbf{A}, \mathbf{W}, \mathbf M) = {\hat \psi}(\mathbf{Q}		\mathbf{v}, \mathbf{Q}\mathbf{A}\mathbf{Q}^T, \mathbf{Q}\mathbf{W} \mathbf{Q}^T, \langle\mathbf Q \rangle \mathbf{M}) \\	
		\mathbf Q \mathbf{\hat T}(\mathbf v, \mathbf{A}, \mathbf{W}, \mathbf M)\mathbf Q^T = \mathbf{\hat T}(\mathbf{Q}		\mathbf{v}, \mathbf{Q}\mathbf{A}\mathbf{Q}^T, \mathbf{Q}\mathbf{W} \mathbf{Q}^T, \langle\mathbf Q \rangle\mathbf{M}) 
	\end{array} ;\forall \mathbf Q \in \mathcal O(2)  
\end{align}
The isotropic scalar- and tensor-valued functions $\hat{\psi}(\mathbf v, \mathbf{A}, \mathbf{W}, \mathbf M)$  and  $\mathbf{\hat T}(\mathbf v, \mathbf{A}, \mathbf{W}, \mathbf M)$ can be represented using Eqs. (\ref{eq 2.2}) and (\ref{eq 2.3}) directly. This theory laid the foundation of constitutive modeling for anisotropic materials. 

Using the aforementioned definition of structural tensors in (\ref{eq 2.6}), Zheng\cite{zheng_two-dimensional_1997} proposed structural tensors for crystallographic and compact point groups in 2D as follows. For an orthonormal vector basis $\{ \mathbf a_1,\mathbf a_2 \}$, two symmetric tensors are defined as 
\begin{align}\label{eq 2.8}
	\begin{array} {ll}
		\mathbf F = \mathbf a_1 \otimes \mathbf a_1 - \mathbf a_2 \otimes \mathbf a_2, & \mathbf H = \mathbf a_1 \otimes \mathbf a_2 + \mathbf a_2 \otimes \mathbf a_1
	\end{array}
\end{align}
and an $n$-th order structural tensor is defined as 
\begin{align} \label{2.9}
	\mathbf P_n = Re(\mathbf a_1+i \mathbf a_2)^{n}= 
	\begin{cases}
		Re(\mathbf F+i \mathbf H)^m, & \text{if } n=2m \\
		Re(\mathbf a_1+i \mathbf a_2) \otimes (\mathbf F+i \mathbf H)^m, & \text{if } n=2m+1
	\end{cases}
\end{align}
where $i=\sqrt{-1}$ is the unit imaginary number and $Re$ is the real part. Therefore, considering $\mathbf a_1 = \mathbf i$ and $\mathbf a_2 = \mathbf j$ leads to $\mathbf P_1 =\mathbf i$, $\mathbf P_2 =\mathbf F$, $\mathbb P_3 =\mathbf i \otimes \mathbf F -\mathbf j \otimes \mathbf H$, $\mathbb P_4 =\mathbf F \otimes \mathbf F -\mathbf H \otimes \mathbf H$, $\mathbb P_6 =\mathbf F \otimes \mathbf F \otimes \mathbf F -(\mathbf F \otimes \mathbf H \otimes \mathbf H +\mathbf H \otimes \mathbf F \otimes \mathbf H +\mathbf H \otimes \mathbf H \otimes \mathbf F)$. The structural tensors provided by Zheng for all crystallographic and compact continuous subgroups of $\mathcal O(2)$ are presented in Table \ref{Table 4}.

Following Boehler and Liu's theory, for the 6 point groups ($\mathcal{C}_1$, $\mathcal{C}_{1v}$, $\mathcal{C}_2$, $\mathcal{C}_{2v}$, $\mathcal{C}_\infty$, $\mathcal{C}_{\infty v}$) with lower-order structural tensors, the representations of a scalar function $\psi(\mathbf C)$ and a tensor function $\mathbf{T}(\mathbf{C})$ with both $\mathbf T$ and $\mathbf C$ as 2nd-order symmetric tensors are presented in Table \ref{Table 3}. However, for the other 6 point groups ($\mathcal{C}_4$, $\mathcal{C}_{4v}$, $\mathcal{C}_3$, $\mathcal{C}_{3v}$, $\mathcal{C}_6$, $\mathcal{C}_{6v}$) with higher-order structural tensors, it is very challenging to find the representation using Boehler and Liu's theory, and this work aims to solve this issue built on the reformulation of Man and Goddard to be introduced next. 

\begin{table}[!h]
	\caption{ Structural tensors proposed by Zheng \cite{zheng_theory_1994} and the authors for some 2D point groups}
	\label{Table 4}
	\centering
	\begin{tabular}{lll}
		\hline
		\centering
		Crystal system (point group) & Zheng's structural tensors \cite{zheng_theory_1994} & Proposed structural tensor set \\ [1.0ex] \hline
		Oblique  $(\mathcal{C}_1)  $  & $\mathbf i$, $\mathbf j$ & $\mathbf i$, $\mathbf j$ \\
		Oblique $\mathcal{(C}_{2})$   & $\mathbf{P}_2, \bm{\varepsilon}$\footnotemark[1] & $\mathbf{P}_2, \bm{\varepsilon}$ \\ [1ex] 
		Rectangular $(\mathcal{C}_{1v})$  & $\mathbf i $ & $\mathbf i $\\
		Rectangular $(\mathcal{C}_{2v})$    &$\mathbf{P}_2 $ &$\mathbf{P}_2 $\\
		[1ex] 
		Square ($\mathcal{C}_{4}  $) & $\mathbb{P}_4, \bm{\varepsilon}$ & $\mathbf i \otimes \mathbf i, \mathbf j \otimes \mathbf j, \bm{\varepsilon}$  \\
		Square ($\mathcal{C}_{4v}  $) & $\mathbb{P}_4$ & $\mathbf i \otimes \mathbf i, \mathbf j \otimes \mathbf j$  \\ [1ex] 
		Hexagonal ($ \mathcal{C}_3   $) & $\mathbb{P}_3, \bm{\varepsilon}$ & $\mathbf v_1,\mathbf v_2,\mathbf v_3, \bm{\varepsilon} $  \\
		Hexagonal ($ \mathcal{C}_{3v}$) & $\mathbb{P}_3$  & $\mathbf v_1,\mathbf v_2,\mathbf v_3 $  \\
		Hexagonal ($ \mathcal{C}_6   $) & $\mathbb{P}_6, \bm{\varepsilon}$ & $\mathbf{M}_1, \mathbf{M}_2, \mathbf{M}_3, \bm{\varepsilon} $  \\
		Hexagonal ($ \mathcal{C}_{6v}   $) & $\mathbb{P}_6$  & $\mathbf{M}_1, \mathbf{M}_2, \mathbf{M}_3 $ \\
		[1ex] 
		Circular ($ \mathcal{C}_{\infty}$)  &  $\bm{\varepsilon}$ &  $\bm{\varepsilon}$ \\
		Circular ($ \mathcal{C}_{\infty v}$)   & $\mathbf{I}$ & $\mathbf{I}$ \\ [1.2ex]
		\hline
	\end{tabular}
	\footnotetext[1]{$\bm{\varepsilon}$ is the second-order permutation tensor in 2D space. See Eq. (\ref{eq 4.1}).}
	\vspace*{-4pt}
\end{table} 

\begin{table}[!h]
	\caption{Structural tensors, functional basis, and tensor generators proposed by Zheng \cite{zheng_theory_1994} for 2nd-order symmetric tensor-valued functions ($\mathcal{C}_1$, $\mathcal{C}_{1v}$, $\mathcal{C}_2$, $\mathcal{C}_{2v}$, $\mathcal{C}_\infty$, $\mathcal{C}_{\infty v}$) }
	\label{Table 3}
	\begin{tabular}{llll}
		\hline
		\centering
		point group   & Structural tensors & Functional basis & Tensor generators \\ [1.0ex] \hline 
		$\mathcal{C}_1  $  & $\mathbf i$, $\mathbf j$ & $\mathbf i.\mathbf C \mathbf i $, $\mathbf j.\mathbf C \mathbf j$, $\mathbf i.\mathbf C \mathbf j $ & $\mathbf i \otimes \mathbf i$, $\mathbf j \otimes \mathbf j$, $\mathbf i \otimes \mathbf j + \mathbf j \otimes \mathbf i$\\
		$\mathcal{C}_{1v}$  & $\mathbf i $ & $tr\mathbf C$, $tr\mathbf C^2$, $\mathbf i.\mathbf C \mathbf i $ & $\mathbf I$, $\mathbf i \otimes \mathbf i$, $\mathbf C$  \\
		$\mathcal{C}_{2}$   & $\mathbf{P}_2, \bm{\varepsilon}$ & $tr\mathbf C$, $tr\mathbf{P}_2 \mathbf C$, $tr\mathbf{P}_2 \mathbf C \bm{\varepsilon}$ & $\mathbf I$, $\mathbf P_2$, $\mathbf P_2 \bm{\varepsilon} $ \\ 
		$\mathcal{C}_{2v}$    &$\mathbf{P}_2 $ &$tr\mathbf C$, $tr\mathbf C^2$,$ tr\mathbf P_2\mathbf C$ &$\mathbf I$, $\mathbf P_2$, $\mathbf C$  \\
		$ \mathcal{C}_{\infty}$  &  $\bm{\varepsilon}$ &$tr\mathbf C$, $tr\mathbf C^2$ & $\mathbf I$, $\mathbf C$, $\mathbf C \bm{\varepsilon} - \bm{\varepsilon} \mathbf C $\\
		$ \mathcal{C}_{\infty v}$   & $\mathbf{I}$ &$tr\mathbf C$, $tr\mathbf C^2$ & $\mathbf I$, $\mathbf C$ \\ [1.2ex]
		\hline
	\end{tabular}
	\vspace*{-4pt}
\end{table}

\subsection{Anisotropic tensor functions: reformulation by Man and Goddard} \label{sec:Man}

In 2018, Man and Goddard \cite{man2018remarks} proposed a new approach and reformulated the representation of anisotropic tensor functions. In their reformulation, the structural tensors do not need to be invariant under symmetry operation of the point group like in Eq. (\ref{eq 2.6}). Instead, the structural tensors can either remain invariant or permute to other structural tensors. This modification allows introducing multiple lower-order (i.e., 1st and 2nd order) structural tensors instead of higher-order structural tensors for many point groups. In the reformulation, the structural tensors only partially preserve the material symmetry so additional symmetry constraints are required for the isotropic extension of the anisotropic tensor functions. 

The reformulated representation theory is briefly summarized below. Note that we use different symbols compared to \cite{man2018remarks}. In this study we present also a reformulation for the scalar-valued tensor functions considering its wide application for constitutive modeling. Consider the scalar- and tensor-valued functions $\psi(\mathbf v, \mathbf{A},\mathbf{W})$ and $\mathbf T (\mathbf v, \mathbf{A},\mathbf{W})$ in Eq. (\ref{eq 2.4}) again. There are three steps to apply the reformulated representation theory:
\begin{enumerate}
	\item Select suitable lower-order structural tensors $\mathbf M$, either vectors or 2nd order tensors, for a specific point group $\mathcal G$. These structural tensors must be either invariant or permute to others under symmmetry operations of $\mathcal G$.
	
	\item Define an isotropic extension of these two anisotropic tensor functions, i.e. ${\psi}(\mathbf v, \mathbf{A},\mathbf{W}) = {\hat \psi}(\mathbf v, \mathbf{A},\mathbf{W},\mathbf{M})$ and $\mathbf {T} (\mathbf v, \mathbf{A},\mathbf{W}) = \mathbf {\hat T} (\mathbf v, \mathbf{A},\mathbf{W},\mathbf{M})$. The representation of ${\hat \psi}$ and $\mathbf {\hat T}$ follows Eqs. (\ref{eq 2.2}) and (\ref{eq 2.3}).
	
	\item Because the structural tensors only partially preserve the material symmetry, additional symmetry constraints are required for ${\hat \psi}$ and $\mathbf {\hat T}$, as
	\begin{align} \label{eq_constraint}
		\begin{array} {ll}
			{\hat \psi}(\mathbf v, \mathbf{A}, \mathbf{W}, \mathbf M) = {\hat \psi}(\mathbf{v}, \mathbf{A}, \mathbf{W}, \langle\mathbf Q \rangle \mathbf{M}) \\	
			\mathbf{\hat T}(\mathbf v, \mathbf{A}, \mathbf{W}, \mathbf M) = \mathbf{\hat T}(\mathbf{v}, \mathbf{A}, \mathbf{W}, \langle\mathbf Q \rangle\mathbf{M}) 
		\end{array} ;\forall \mathbf Q \in \mathcal G  
	\end{align}
	In fact, only the point group generators $\mathcal G^*$ need to be imposed in Eq. (\ref{eq_constraint}) instead of the whole point group $\mathcal G$. A proof is provided in the SI. This would greatly simplify the constraint equations. The point group generators $\mathcal G^*$ of all 2D point groups are listed in Table S11 of the Supplementary Information.      
\end{enumerate}

In short, the core idea is to first identify isotropic tensor functions $ \hat{\psi}$  and $\hat{\mathbf{T}}$ in an augmented function space by initially relaxing the symmetry constraints on the structural tensors. The function space is then progressively narrowed by reintroducing the appropriate symmetry constraints using Eq.~(\ref{eq_constraint}). While determining the specific functional forms remains challenging, the reformulation proposed by Man and Goddard provides a promising path forward.

Man and Goddard \cite{man2018remarks} only introduced the mathematical framework and gave a few examples. In order to apply this theory for constitutive modeling, suitable lower-order structural tensors remain to be discovered for each point group. In this work, we aim to fill this knowledge gap for all 2D point groups and derive explicit formulae for the representation of scalar- and 2nd-order symmetric tensor-valued anisotropic functions. 

\section{Lower-order structural tensor set} \label{sec:struct_tensor_set}

This section introduces a method to discover lower-order structural tensors for all 2D point groups that are suitable for Man and Goddard's reformulation. 

Provided that these lower-order structural tensors are either invariant or permute to others under any symmetry operation in $\mathcal G$, we introduce a new concept \textit{structural tensor set} $\{ \mathbf{M}_i \}$, which has the following mathematical definition 
\begin{align} \label{3.1}
	\mathcal G_s=\{ \mathbf Q \in \mathcal O (2)\ |\  \{\langle \mathbf Q\rangle \mathbf M_i \} = \{\mathbf M_i \} ;\ \ i=1,2,...,s\}
\end{align}
where $\mathcal G_s$ is a subgroup of $\mathcal G$, i.e. $\mathcal{G}_s \leq
 \mathcal{G}$. Herein, we say that the structural tensor set $\{ \mathbf M_i \}$ characterizes the group $\mathcal {G} _s$.  The traditional definition of structural tensors Eq. (\ref{eq 2.6}) requires each $\mathbf{M}_i$ to be invariant $\forall Q \in \mathcal{G}$. In contrast, Eq. (\ref{3.1}) only requires the whole structural tensor set to be invariant $\forall Q \in \mathcal{G}_s$. Note that herein we do not mandate $\mathcal{G}_s = \mathcal{G}$, although it is still strongly recommended (see the Remark 3.2 below). 

\begin{figure}[!h]
	\centering\includegraphics[width=4 in]{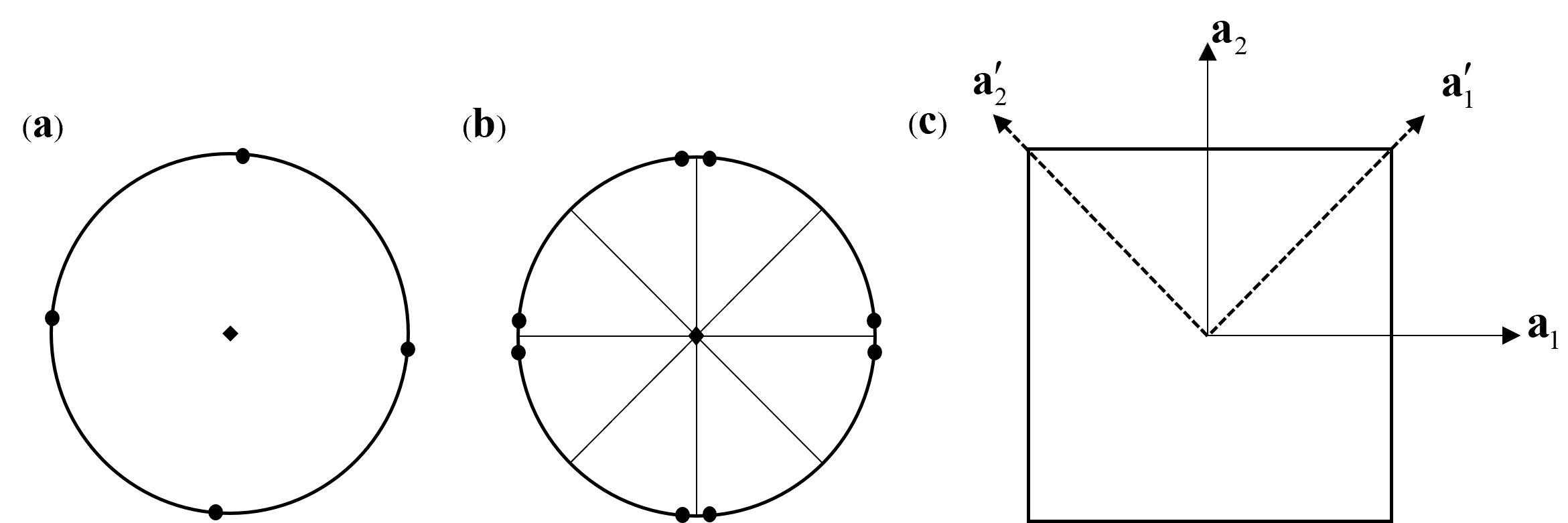}
	\centering \caption{Graphical representation of the 2D point group (a) $\mathcal C_4$ and (b) $\mathcal C_{4v}$. (c) Illustration of various vectors used to define structural tensor set for $\mathcal C_4$ and $\mathcal C_{4v}$.}
	\label{figure 1}
\end{figure}

Certainly, the structural tensor sets are non-unique and one has to discover a convenient set for each point group. It usually takes laborious work and a lengthy trial-and-error process to find a complete structural tensor set. Generally, one can choose typical high symmetry directions, lines, and planes to construct the lower-order structural tensors. Take the 2D point groups  $\mathcal C_4$ and $\mathcal C_{4v}$ in Figure \ref{figure 1} as an example. 
Considering two orthonormal vectors $\{\mathbf a_1, \mathbf a_2\} = \{\mathbf i, \mathbf j\}$, one option of the structural tensor set for $\mathcal C_{4v}$ and $\mathcal C_{4}$ are $\{\mathbf i \otimes \mathbf i,\mathbf j \otimes \mathbf j\}$ and $\{\mathbf i \otimes \mathbf i,\mathbf j \otimes \mathbf j, \bm{\varepsilon} \}$, respectively. If we adopt two different orthonormal vectors $\{{\mathbf a}'_1, {\mathbf a}'_2\} = \{(\mathbf{i}+\mathbf{j})/ \sqrt 2, (\mathbf{j}-\mathbf{i})/ \sqrt 2 \}$ in Figure \ref{figure 1}(c), another option of the structural tensor set would be $\{\mathbf{a}'_1 \otimes \mathbf{a}'_1,\mathbf{a}'_2 \otimes \mathbf{a}'_2 \}$ with  for point group $\mathcal C_{4v}$ and $\{\mathbf{a}'_1 \otimes \mathbf{a}'_1,\mathbf{a}'_2 \otimes \mathbf{a}'_2, \bm{\varepsilon} \}$ for point group $\mathcal C_{4}$. Note that adding a 2D permutation tensor $\bm{\varepsilon}$ to the structural tensor set would break the reflection symmetry; meanwhile, $\mathcal C_{4v}$ degenerates to $\mathcal C_{4}$. 

\textbf{Remark 3.1.} In Eq. (\ref{3.1}), $\mathcal G_s$ can be a proper subgroup of $\mathcal G$, but not a supergroup. In other words, $\mathcal G_s$ can have lower symmetry than $\mathcal G$, but not higher. If $\mathcal G_s  < \mathcal G$, additional symmetry constraints must be imposed through Eq. (\ref{eq_constraint}). However, if $\mathcal G_s > \mathcal G$, the structural tensor set introduces hidden (or redundant) symmetry to the representation. For example, if we first define the structural tensor set for $\mathcal C_{6v}$ and use the same set for $\mathcal C_{3v}$, a hidden inversion symmetry may be added to the representation of $\mathcal C_{3v}$. In contrast, one can always start with the representation of the subgroup (e.g., $\mathcal C_{3v}$ here), and adds more symmetry constraints in Eq. (\ref{eq_constraint}) to obtain the representation of $\mathcal C_{6v}$.

\textbf{Remark 3.2.} For the convenience of modeling, we strongly suggest choosing a structural tensor set that characterizes the point group, i.e., $\mathcal{G}_s = \mathcal{G}$. This would ease the burden to impose symmetry constraints  (\ref{eq_constraint}). In other words, if we over-relax the symmetry constraints of structural tensors at beginning, the augmented function space $\hat \psi$ and $\hat{\mathbf T}$ will be very large, and then it would take more efforts to narrow down the function space by adding additional constraints (\ref{eq_constraint}).     

\textbf{Remark 3.3.} The structural tensor set $\{\mathbf{M}_i \}$ must be complete and invariant $\forall \mathbf{Q} \in \mathcal{G}_s$. This necessary condition needs to be verified for any newly proposed structural tensor sets. In addition, this condition is useful when constructing a structural tensor set. For example, starting with one structural tensor $\mathbf{M}$, one can find a few other structural tensors by performing symmetry transformation $ \langle \mathbf Q\rangle \mathbf M  $ for $\forall \mathbf{Q} \in \mathcal{G}_s$. 

In brief summary, for the six point groups ($\mathcal{C}_1$, $\mathcal{C}_{1v}$, $\mathcal{C}_2$, $\mathcal{C}_{2v}$, $\mathcal{C}_\infty$, $\mathcal{C}_{\infty v}$) with lower-order structural tensors, the definition in Eq. (\ref{3.1}) actually degenerates to the traditional structural tensor definition Eq. (\ref{eq 2.6}) when $\mathcal{G} = \mathcal{G}_s$. Hence for these six groups, we simply adopt Zheng's structural tensors in Table \ref{Table 4}. In these cases, Boehler and Liu's formulation in Eqs. (\ref{eq 2.2}) and (\ref{eq 2.3}) should be used to derive the representations with the required functional basis and tensor generators given in Table \ref{Table 3}.  In contrast, for the six groups ($\mathcal{C}_4$, $\mathcal{C}_{4v}$, $\mathcal{C}_3$, $\mathcal{C}_{3v}$, $\mathcal{C}_6$, $\mathcal{C}_{6v}$) with higher-order structural tensors, we propose lower-order structural tensor sets instead of Zheng's structural tensors, as tabulated in Table \ref{Table 4}. For these six groups, the reformulation of Man and Goddard should be used. The representations of scalar- and tensor-valued functions for these six 2D point groups are presented in Sections \ref{sec:C4}, \ref{sec:C3}, and \ref{sec:C6} below.

\section{Representation of tensor functions for $\mathcal C_4$ and $\mathcal C_{4v}$} \label{sec:C4}
\subsection{Group $\mathcal C_4$}
We will follow the three steps outlined in Section \ref{sec:Man} to derive the representations using the reformulation of Man and Goddard.

For the point group $\mathcal C_4$ illustrated in Figure \ref{figure 1}(a), we propose a structural tensor set $\{ \mathbf{M}_1, \mathbf{M}_2, \bm{\varepsilon} \}$ comprising of three tensors, as
\begin{align} \label{eq 4.1}
	\begin{array}{l}
		{\mathbf{M}_1} = \mathbf i \otimes \mathbf i = \left[ {\begin{array}{*{20}{c}}
				1&0\\
				0&0
		\end{array}} \right], \ 
		{\mathbf{M}_2} = \mathbf j \otimes \mathbf j = \left[ {\begin{array}{*{20}{c}}
				0&0\\
				0&1
		\end{array}} \right],\ 
		\bm{\varepsilon}  = \left[ {\begin{array}{*{20}{c}}
				0&1\\
				{ - 1}&0
		\end{array}} \right]
	\end{array}   
\end{align}
where $\mathbf M_1 = \mathbf i \otimes  \mathbf i$ and $\mathbf M_2 = \mathbf j \otimes  \mathbf j$ are defined based on the two orthonormal vector bases and $\bm{\varepsilon}$ is introduced to eliminate the reflection symmetry. 

The representation of a 2nd-order symmetric-tensor-valued function $\mathbf{T(C)}$ is derived first, where the argument $ \mathbf C$ is also a 2nd-order symmetric tensor (e.g., the right Cauchy-Green tensor in continuum mechanics). We first define an isotropic extension $\mathbf{\hat T}$ for $\mathbf T$ using the structural tensor set, as 
\begin{align} \label{eq 4.2}
	\mathbf{T(C) = \hat T(C},{\mathbf{M}_1},{\mathbf{M}_2},\bm{\varepsilon})   
\end{align}	
Using this point group's generator $\mathcal{G^*} = \{ {\mathbf Q^{\pi /2}}\}$ where $ \mathbf Q^{\pi /2} = \begin{bmatrix}
0 & 1 \\
-1 & 0
\end{bmatrix} $, (see Table S11 in Supplementary Information), we can find that
\begin{align} \label{eq 4.3}
	\begin{array}{l}
		{\mathbf Q^{\pi /2}}{\mathbf{M}_1}{({\mathbf Q^{\pi /2}})^T} = {\mathbf{M}_2}\\
		{\mathbf Q^{\pi /2}}{\mathbf{M}_2}{({\mathbf Q^{\pi /2}})^T} = {\mathbf{M}_1}\\
		{\mathbf Q^{\pi /2}}{\bm{\varepsilon}}{({\mathbf Q^{\pi /2}})^T} = {\bm{\varepsilon}}
	\end{array}	   
\end{align}
Obviously, the two structural tensors $\mathbf{M}_1$ and $\mathbf M_2$ are not invariant under $\mathbf Q^{\pi /2}$ but permute to each other. This is the case when additional symmetry constraints need to be imposed through Eq. (\ref{eq_constraint}) following the reformulation of Man and Goddard. 

In order to find the representation of $\mathbf T(\mathbf C)= \mathbf{\hat T(C},{\mathbf{M}_1},{\mathbf{M}_2},\bm{\varepsilon})$, considering that $\mathbf C, {\mathbf{M}_1}, {\mathbf{M}_2}$ are $2^{nd}$-order symmetric tensors and $\bm{\varepsilon}$ is a skew-symmetric tensor and using Tables \ref{Table 1} and \ref{Table 2}, we can obtain
\begin{align}\label{eq 4.4}
	\begin{aligned}& \mathbf{T(C)} =\mathbf{ \hat T(C},{\mathbf{M}_1},{\mathbf{M}_2},\bm{\varepsilon}) = {\alpha _0}\mathbf I + {\alpha _1}\mathbf C + {\alpha _2}{\mathbf{M}_1} + {\alpha _3}{\mathbf{M}_2} + {\alpha _4}(\mathbf C \bm{\varepsilon} - \bm{\varepsilon} \mathbf C)\\ &
		\quad + {\alpha _5}({\mathbf{M}_1}{\bm{\varepsilon}} - {\bm{\varepsilon}}{\mathbf{M}_1}) + {\alpha _6}({\mathbf{M}_2}{\bm{\varepsilon}} - {\bm{\varepsilon}}{\mathbf{M}_2}) 
	\end{aligned}
\end{align}
where $\alpha_i$ is a scalar coefficient function given by  
\begin{align}\label{eq 4.5}
	\begin{aligned}& {\alpha _i} = {\alpha _i}(tr \mathbf C,tr{ \mathbf C^2},tr{\mathbf{M}_1},tr{\mathbf{M}_1^2},tr{\mathbf{M}_2},tr{\mathbf{M}_2^2},tr{\bm{\varepsilon}^2},tr(\mathbf C{\mathbf{M}_1}),tr(\mathbf C{\mathbf{M}_2}),\\&
		\quad tr({\mathbf{M}_1}{\mathbf{M}_2}),tr(\mathbf C{\mathbf{M}_1}{\bm{\varepsilon}}),tr(\mathbf C{\mathbf{M}_2}{\bm{\varepsilon}}),tr({\mathbf{M}_1}{\mathbf{M}_2}{\bm{\varepsilon}})) 
	\end{aligned}
\end{align}
Note that there are redundant terms in (\ref{eq 4.4}) and (\ref{eq 4.5}). After simplification, the representation is re-written as 
\begin{align}\label{eq 4.6}
	\begin{aligned} &\mathbf{T(C) = \hat T(C},{\mathbf{M}_1},{\mathbf{M}_2},{\bm{\varepsilon}}) = {\alpha _0}\mathbf C + {\alpha _1}{\mathbf{M}_1} + {\alpha _2}{\mathbf{M}_2} + {\alpha _3}(\mathbf C{\bm{\varepsilon}} - {\bm{\varepsilon}}\mathbf C) \\ & \quad \quad \ \  + {\alpha _4}({\mathbf{M}_1}{\bm \varepsilon} - {\bm \varepsilon}{\mathbf{M}_1})   
	\end{aligned}
\end{align}
and 
\begin{align} \label{eq 4.7}
	\begin{aligned} & {\alpha _i} = {\alpha _i}(tr(\mathbf C{\mathbf{M}_1}),tr(\mathbf C{\mathbf{M}_2}),tr(\mathbf C{\mathbf{M}_1}{\bm{\varepsilon}}))\\ &
		\quad ={\tilde{\alpha} _i}(\mathbf{C},\mathbf M_1, \mathbf M_2, \bm{\varepsilon});\ i = 0,1,...,4   
	\end{aligned}
\end{align}
Moreover, additional constraints on (\ref{eq 4.6}) are imposed as based on (\ref{eq_constraint}), as 
\begin{align} \label{eq 4.8}
	\mathbf{\hat T(C},{\mathbf{M}_1},{\mathbf{M}_2},{\bm{\varepsilon}}) = \mathbf{\hat T(C},{\mathbf Q}{\mathbf{M}_1}{\mathbf Q^T},{\mathbf Q}{\mathbf{M}_2}{\mathbf Q^T},{\mathbf Q}{\bm{\varepsilon}}{\mathbf Q^T}) ; \quad \forall \mathbf Q \in \mathcal G^*  
\end{align}
Hence, for ${\mathbf Q} = {\mathbf Q^{\pi /2}}$, the constraint is
\begin{align} \label{eq 4.9}
	\mathbf{\hat T(C},{\mathbf{M}_1},{\mathbf{M}_2},{\bm{\varepsilon}}) = \mathbf{\hat T(C},{\mathbf{M}_2},{\mathbf{M}_1},{\bm{\varepsilon}})  
\end{align}
Finally, by substituting (\ref{eq 4.6}) in (\ref{eq 4.9}), we can obtain the constraints on the scalar coefficient functions as follows
\begin{align} \label{eq 4.10}
	\begin{aligned} &  \tilde{\alpha} _0(\mathbf C,{\mathbf{M}_1},{\mathbf{M}_2},{\bm{\varepsilon}}) = \tilde{\alpha} _0(\mathbf C,{\mathbf{M}_2},{\mathbf{M}_1},{\bm{\varepsilon}}){\rm{ }}\\&
		\tilde{\alpha} _1(\mathbf C,{\mathbf{M}_1},{\mathbf{M}_2},{\bm{\varepsilon}}) = \tilde{\alpha} _2(\mathbf C,{\mathbf{M}_2},{\mathbf{M}_1},{\bm{\varepsilon}})\\&
		\tilde{\alpha} _3(\mathbf C,{\mathbf{M}_1},{\mathbf{M}_2},{\bm{\varepsilon}}) = \tilde{\alpha} _3(\mathbf C,{\mathbf{M}_2},{\mathbf{M}_1},{\bm{\varepsilon}})\\&
		\tilde{\alpha} _4(\mathbf C,{\mathbf{M}_1},{\mathbf{M}_2},{\bm{\varepsilon}}) =- \tilde{\alpha} _4(\mathbf C,{\mathbf{M}_2},{\mathbf{M}_1},{\bm{\varepsilon}})
	\end{aligned}  
\end{align}

The representation of a scalar function $\psi(\mathbf C)$ can be derived following a similar procedure. We first define an isotropic extension as $\psi(\mathbf C) = \hat \psi(\mathbf C, \mathbf M_1, \mathbf M_2, \bm{\varepsilon})$. The representation of $\hat \psi$ follows the same form of (\ref{eq 4.7}), as
\begin{align}\label{eq psi_C4}
	\begin{array}{l}
		\psi(\mathbf C) = \hat \psi(\mathbf C, \mathbf M_1, \mathbf M_2, \bm{\varepsilon})
		= \hat \psi(tr(\mathbf C\mathbf M_1),tr(\mathbf C \mathbf M_2),tr(\mathbf C \mathbf M_1 \bm{\varepsilon}))
	\end{array}
\end{align}
Next, similar to (\ref{eq 4.8}), additional constraints must be imposed as
\begin{align} \label{eq 4.8 psi}
	\mathbf{\hat \psi(C},{\mathbf{M}_1},{\mathbf{M}_2},{\bm{\varepsilon}}) = \mathbf{\hat \psi(C},{\mathbf Q}{\mathbf{M}_1}{\mathbf Q^T},{\mathbf Q}{\mathbf{M}_2}{\mathbf Q^T},{\mathbf Q}{\bm{\varepsilon}}{\mathbf Q^T}) ; \quad \forall \mathbf Q \in \mathcal G^*  
\end{align}
For ${\mathbf Q} = {\mathbf Q^{\pi /2}}$, this constraint turns out to be $ \mathbf{\hat \psi(C},{\mathbf{M}_1},{\mathbf{M}_2},{\bm{\varepsilon}}) = \mathbf{\hat \psi(C},{\mathbf{M}_2},{\mathbf{M}_1},{\bm{\varepsilon}}) $, or more explicitly,
\begin{align}\label{eq psi_constraint_C4}
	\begin{array}{l}
		\hat \psi(tr(\mathbf C\mathbf M_1),tr(\mathbf C \mathbf M_2),tr(\mathbf C \mathbf M_1 \bm{\varepsilon}))
		= \hat \psi(tr(\mathbf C\mathbf M_2),tr(\mathbf C \mathbf M_1),-tr(\mathbf C \mathbf M_1 \bm{\varepsilon}))
	\end{array}
\end{align}

\textbf{Remark 4.1.} One may take advantage of identities to reduce the number of functional bases and tensor generators to simplify the representations (not necessarily irreducible). For example, the last two terms in (\ref{eq 4.4}) have the relation $(\mathbf{M}_1 \bm{\varepsilon}-\bm{\varepsilon}\mathbf{M}_1)=-(\mathbf{M}_2 \bm{\varepsilon}-\bm{\varepsilon}\mathbf{M}_2)$. Hence, this would eliminate one of them in (\ref{eq 4.6}) and also give rise to the negative sign in the last equation of (\ref{eq 4.10}). In addition, take the identity $tr(\mathbf C \mathbf M_1 \bm{\varepsilon})= -tr(\mathbf C \mathbf M_2 \bm{\varepsilon})$ in (\ref{eq 4.5}) as another example. Clearly, $tr(\mathbf C \mathbf M_2 \bm{\varepsilon})$ is a redundant term in the functional bases. Consequently, this also replaces the argument $tr(\mathbf C \mathbf M_2 \bm{\varepsilon})$ with $-tr(\mathbf C \mathbf M_1 \bm{\varepsilon})$ on the right hand side of (\ref{eq psi_constraint_C4}). 

\textbf{Remark 4.2.} When one eliminates the redundant terms in (\ref{eq 4.4}) and (\ref{eq 4.5}), there is no unique way to do it. The leftover functional bases and tensor generators are chosen based on users' preference. If a different set of functional bases and tensor generators are chosen, the detailed representations will be different. This rule applies to other point groups as well. 

\textbf{Remark 4.3.} Some constraints on the scalar coefficient functions $\tilde{\alpha}_i$ are redundant and should be removed. The fact is that not all these constraints are independent. For example, in (\ref{eq 4.10}), we have removed a constraint $\tilde{\alpha} _2(\mathbf C,{\mathbf{M}_1},{\mathbf{M}_2},{\bm{\varepsilon}}) = \tilde{\alpha} _1(\mathbf C,{\mathbf{M}_2},{\mathbf{M}_1},{\bm{\varepsilon}})$ because it is equivalent to the 2nd equation in (\ref{eq 4.10}). Other point groups have the same issue and careful attention must be paid to remove redundancy.  

\subsection{Group $\mathcal C_{4v}$}
The representation of tensor functions for $\mathcal C_{4v}$ is simpler than $\mathcal C_{4}$. The derivation procedure is quite similar. In the case of $\mathcal C_{4v}$, a structural tensor set $\{ \mathbf M_1, \mathbf M_2 \}$ is chosen with ${\mathbf{M}_1}$ and ${\mathbf{M}_2}$ defined in (\ref{eq 4.1}). By eliminating all terms of $\bm{\varepsilon}$ in (\ref{eq 4.4}) and (\ref{eq 4.5}) and removing redundant terms, we can obtain the representation of $\mathbf{T(C)}$ as    
\begin{align} \label{eq 4.11}
	\mathbf{T(C) = \hat T(C},{\mathbf{M}_1},{\mathbf{M}_2}) = {\alpha _0}\mathbf C + {\alpha _1}{\mathbf{M}_1} + {\alpha _2}{\mathbf{M}_2}   
\end{align}
and 
\begin{align} \label{eq 4.12}
	{\alpha _i} = {\alpha _i}({{tr}}{\mathbf C^2},{{tr(}}\mathbf C{\mathbf{M}_1}),{{ tr(}}\mathbf C{\mathbf{M}_2}))={\tilde{\alpha}_i}(\mathbf{C},\mathbf M_1, \mathbf M_2), \ i=0,1,2   
\end{align}

According to (\ref{eq_constraint}), additional constraints must be imposed to the coefficient functions. Using the point group's generators $\mathcal G^* = \{ \mathbf Q^{\pi /2},\mathbf m_{10}\}$ (see Table S11 in Supplementary Information), we can find that
\begin{align} \label{eq 4.13}
	\begin{aligned} & {\mathbf Q^{\pi /2}}{\mathbf{M}_1}{({\mathbf Q^{\pi /2}})^T} = {\mathbf{M}_2},\ {\rm{ }}{\mathbf m_{10}}{\mathbf{M}_1}{{\mathbf m_{10}^T}} = {\mathbf{M}_1}\\ &
		{\mathbf Q^{\pi /2}}{\mathbf{M}_2}{({\mathbf Q^{\pi /2}})^T} = {\mathbf{M}_1},\ {\rm{ }}{\mathbf m_{10}}{\mathbf{M}_2}{{\mathbf m_{10}^T}} = {\mathbf{M}_2}
	\end{aligned}
\end{align}
where  $ \mathbf Q^{\pi /2} = \begin{bmatrix}
0 & 1 \\
-1 & 0
\end{bmatrix} $ and $\mathbf m_{10} = \begin{bmatrix}
-1 & 0 \\
0 & 1
\end{bmatrix}$ represent the $\pi/2$ rotation symmetry and vertical reflection symmetry, respectively.  According to (\ref{eq 4.13}), the two structural tensors $\mathbf M_1$ and $\mathbf M_2$ are invariant under reflection $\mathbf m_{10}$ and permuting to each other under rotation $\mathbf Q^{\pi/2}$. Hence, we only need to impose additional $\mathbf Q^{\pi /2}$ symmetry to (\ref{eq 4.11}), i.e., $\mathbf{\hat T(C},{\mathbf{M}_1},{\mathbf{M}_2}) = \mathbf{\hat T(C},{\mathbf{M}_2},{\mathbf{M}_1})$. Further, such constraints require that the coefficient functions to satisfy 
\begin{align}
	\begin{aligned} &{\tilde\alpha _0(\mathbf C,{\mathbf{M}_1},{\mathbf{M}_2}) = \tilde\alpha _0(\mathbf C,{\mathbf{M}_2},{\mathbf{M}_1})} \\ & 
		\tilde\alpha _1(\mathbf C,{\mathbf{M}_1},{\mathbf{M}_2}) = \tilde\alpha _2(\mathbf C,{\mathbf{M}_2},{\mathbf{M}_1})    
	\end{aligned}
\end{align}

The representation of scalar function is presented in what follows. Similar to (\ref{eq 4.12}), the representation of a scalar function $\psi(\mathbf C)$ is 
\begin{align}\label{eq psi_C4v}
	\begin{array}{l}
		\psi(\mathbf C) = \hat \psi(\mathbf C, \mathbf M_1, \mathbf M_2)
		= \hat \psi(tr\mathbf C^2, tr(\mathbf C\mathbf M_1),tr(\mathbf C \mathbf M_2))
	\end{array}
\end{align}
Moreover, the additional constraint $\hat \psi(\mathbf C, \mathbf M_1, \mathbf M_2) = \hat \psi(\mathbf C, \mathbf M_2, \mathbf M_1) $ requires that
\begin{align}\label{eq psi_constraint_C4v}
	\begin{array}{l}
		\hat \psi(tr\mathbf C^2, tr(\mathbf C\mathbf M_1),tr(\mathbf C \mathbf M_2))
		= \hat \psi(tr\mathbf C^2, tr(\mathbf C\mathbf M_2),tr(\mathbf C \mathbf M_1))
	\end{array}
\end{align}

\section{Representation of tensor functions for $\mathcal C_3$ and $\mathcal C_{3v}$} \label{sec:C3}
\subsection{Group $\mathcal C_3$}
\begin{figure}[!h]
	\centering\includegraphics[width=4 in]{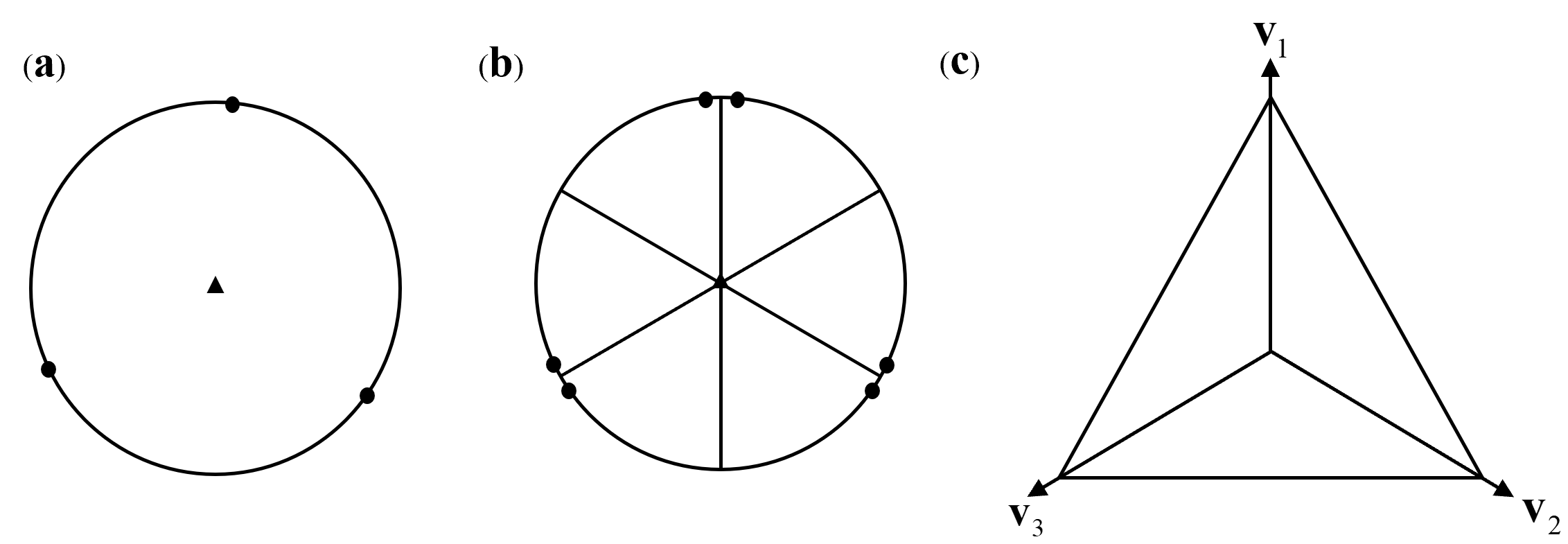}
	\centering \caption{Graphical representation of the 2D point group (a) $\mathcal C_3$ and (b) $\mathcal C_{3v}$. (c) Illustration of vectors used to define structural tensor set for point groups $\mathcal C_3$ and $\mathcal C_{3v}$.} 
	\label{figure 2}
\end{figure}

As illustrated in Figure \ref{figure 2}, the point group $\mathcal C_3$ has a 3-fold rotation symmetry and no reflection symmetry. In order to construct its structural tensor set, we introduce three vectors $\mathbf v_1$, $\mathbf v_2$, and $\mathbf v_3$ in Figure \ref{figure 2}c to characterize the symmetry. Mathematically, they are expressed as
\begin{align} \label{eq structural_C3}
	\begin{aligned} & \mathbf{v}_1 = \begin{bmatrix}
			0\\
			1
		\end{bmatrix}, \mathbf{v}_2 = \begin{bmatrix}
			\sqrt{3}/2\\
			-1/2
		\end{bmatrix}, \mathbf{v}_3 = \begin{bmatrix}
			-\sqrt{3}/2\\
			-1/2
		\end{bmatrix}
	\end{aligned}
\end{align}
The whole structural tensor set for $\mathcal C_3$ is defined as $\{ \mathbf{v}_1, \mathbf{v}_2, \mathbf{v}_3, \bm{\varepsilon} \}$, where $\bm{\varepsilon}$ is introduced to eliminate the reflection symmetry. 

The point group generator of $\mathcal C_3$ is $\mathcal{G}^* = \{\mathbf{Q}^{2\pi /3}\}$, where $\mathbf{Q}^{2\pi /3} = \begin{bmatrix}
-1/2 & \sqrt{3}/2 \\
-\sqrt{3}/2 & -1/2
\end{bmatrix}$ denotes $2 \pi/3$ rotation symmetry (see Table S11 in Supplementary Information). Under $\mathbf{Q}^{2\pi /3}$ transformation, the three structural vectors $\mathbf v_1$, $\mathbf v_2$, and $\mathbf v_3$ permute to each other, whereas $ \bm{\varepsilon}$ is invariant, i.e.,     
\begin{align} \label{eq g*_C3}
	\begin{aligned} & \mathbf{Q}^{2\pi /3}\mathbf v_1 =\mathbf v_2, \  
		\mathbf{Q}^{2\pi /3}\mathbf v_2 =\mathbf v_3, \ 
		\mathbf{Q}^{2\pi /3}\mathbf v_3 =\mathbf v_1, \\ & 
		\mathbf{Q}^{2\pi /3}\bm{\varepsilon}(\mathbf{Q}^{2\pi /3})^T = \bm{\varepsilon}
	\end{aligned}
\end{align}

The representation of a tensor function $\mathbf T(\mathbf C)$ is presented. Firstly, we define its isotropic extension $\mathbf{\hat T}$ using the structural tensor set, as 
\begin{align} \label{eq 5.2}
	\mathbf{T(C)} = \mathbf{\hat T}(\mathbf{C},{\mathbf{v}_1},{\mathbf{v}_2},{\mathbf{v}_3}, \bm{\varepsilon})    
\end{align}
Then, using Tables \ref{Table 1} and \ref{Table 2}, the representation of $\mathbf T(\mathbf C)= \mathbf{\hat T(C},{\mathbf{v}_1},{\mathbf{v}_2},{\mathbf{v}_3}, \bm{\varepsilon} )$ is obtained as 	
\begin{align} \label{eq generators_redundant_C3}
	\begin{aligned} &\mathbf{T(C)} = \mathbf{\hat T}(\mathbf{C},{\mathbf{v}_1},{\mathbf{v}_2},{\mathbf{v}_3}, \bm{\varepsilon}) = \alpha_0\mathbf{I} + \alpha _1\mathbf{C} + \alpha_2 (\mathbf{v}_1 \otimes \mathbf{v}_1)+ \alpha_3(\mathbf{v}_2 \otimes \mathbf{v}_2)\\ & + \alpha _4(\mathbf{v}_3 \otimes \mathbf{v}_3)+  \alpha_5(\mathbf{v}_1 \otimes \mathbf{v}_2 + \mathbf{v}_2 \otimes \mathbf{v}_1) + \alpha_6(\mathbf{v}_1 \otimes \mathbf{v}_3 + \mathbf{v}_3 \otimes \mathbf{v}_1)+ \\ & \alpha_7(\mathbf{v}_2 \otimes \mathbf{v}_3 + \mathbf{v}_3 \otimes \mathbf{v}_2)  + \alpha_8(\mathbf{C}\bm{\varepsilon}-\bm{\varepsilon} \mathbf{C} )  + \alpha_9 (\mathbf{v}_1 \otimes \bm{\varepsilon} \mathbf{v}_1 + \bm{\varepsilon} \mathbf{v}_1 \otimes \mathbf{v}_1 )\\ & + \alpha_{10} (\mathbf{v}_2 \otimes \bm{\varepsilon} \mathbf{v}_2 + \bm{\varepsilon} \mathbf{v}_2 \otimes \mathbf{v}_2 ) + \alpha_{11} (\mathbf{v}_3 \otimes \bm{\varepsilon} \mathbf{v}_3 + \bm{\varepsilon} \mathbf{v}_3 \otimes \mathbf{v}_3 )
	\end{aligned}
\end{align}
where 
\begin{align} \label{eq basis_redundant_C3}
	\begin{aligned} & \alpha _i = \alpha _i(tr\mathbf{C}, tr\mathbf{C}^2, \mathbf{v}_1.\mathbf{C}\mathbf{v}_1, \mathbf{v}_2.\mathbf{C}\mathbf{v}_2, \mathbf{v}_3.\mathbf{C}\mathbf{v}_3, , \mathbf{v}_1.\mathbf{C}\mathbf{v}_2, \mathbf{v}_1.\mathbf{C}\mathbf{v}_3, \mathbf{v}_2.\mathbf{C}\mathbf{v}_3,\\ & \quad \mathbf{v}_1.\mathbf{v}_1, \mathbf{v}_2.\mathbf{v}_2, \mathbf{v}_3.\mathbf{v}_3, \mathbf{v}_1.\mathbf{v}_2, \mathbf{v}_1.\mathbf{v}_3, \mathbf{v}_2.\mathbf{v}_3, tr\bm{\varepsilon}^2, \mathbf{v}_1.\bm{\varepsilon} \mathbf{v}_2, \mathbf{v}_1.\bm{\varepsilon} \mathbf{v}_3, \mathbf{v}_2.\bm{\varepsilon} \mathbf{v}_3, \\ & \quad (\mathbf{C} \mathbf{v}_1).(\bm{\varepsilon} \mathbf{v}_1), ( \mathbf{C} \mathbf{v}_2).(\bm{\varepsilon} \mathbf{v}_2),  (\mathbf{C} \mathbf{v}_3).(\bm{\varepsilon} \mathbf{v}_3) )
	\end{aligned}
\end{align}	
Next, after eliminating the redundant terms, (\ref{eq generators_redundant_C3}) and (\ref{eq basis_redundant_C3}) are further simplified as follows.
\begin{align} \label{eq generators_C3}
	\begin{aligned} & \mathbf{T(C)} = \mathbf{\hat T}(\mathbf{C},{\mathbf{v}_1},{\mathbf{v}_2},{\mathbf{v}_3}, \bm{\varepsilon}) = \alpha_0 \mathbf C + \alpha_1 (\mathbf{v}_1 \otimes \mathbf{v}_1)+ \alpha_2(\mathbf{v}_2 \otimes \mathbf{v}_2)\\ & \quad \quad \ \  + \alpha _3(\mathbf{v}_3 \otimes \mathbf{v}_3)   + \alpha_4(\mathbf{C}\bm{\varepsilon}-\bm{\varepsilon} \mathbf{C} ) 
	\end{aligned}
\end{align}
where 
\begin{align}  \label{eq basis_C3}
	\begin{array}{ll}
		{\alpha _i = \alpha _i( \mathbf v_1.\mathbf{C} \mathbf v_1, \mathbf v_2.\mathbf{C} \mathbf v_2, \mathbf v_3.\mathbf{C} \mathbf v_3)}\\
		\quad \ =\tilde{\alpha}_i(\mathbf C, \mathbf v_1,\mathbf v_2, \mathbf v_3, \bm{\varepsilon}); \quad i=0,1,...,4   
	\end{array}
\end{align}
Finally, we must impose additional constraints (\ref{eq_constraint}) to the representation in (\ref{eq generators_C3}), as 
\begin{align} \label{eq 5.7}
	\mathbf{\hat T(C},{\mathbf{v}_1},{\mathbf{v}_2},{\mathbf{v}_3},\bm{\varepsilon} ) = \mathbf{\hat T(C},{\mathbf Q}{\mathbf{v}_1},{\mathbf Q}{\mathbf{v}_2},{\mathbf Q}{\mathbf{v}_3}, {\mathbf Q}\bm{\varepsilon}{\mathbf Q}^T) ; \quad \forall \mathbf Q \in \mathcal G^*  
\end{align}
When ${\mathbf Q} = {\mathbf Q^{2\pi /3}}$, (\ref{eq 5.7}) leads to
\begin{align} \label{eq 5.8}
	\mathbf{\hat T(C},{\mathbf{v}_1},{\mathbf{v}_2},{\mathbf{v}_3},\bm{\varepsilon}) =  \mathbf{\hat T(C},{\mathbf{v}_2},{\mathbf{v}_3},{\mathbf{v}_1},\bm{\varepsilon} ) 
\end{align}
By substituting (\ref{eq generators_C3}) to (\ref{eq 5.8}), we can obtain the following constraints on the scalar coefficient functions:
\begin{align} \label{eq constraint_C3}
	\begin{aligned} &
		\tilde\alpha_0(\mathbf{C},\mathbf{v}_1,\mathbf{v}_2,\mathbf{v}_3,\bm{\varepsilon}) = \tilde\alpha_0(\mathbf{C},\mathbf{v}_2,\mathbf{v}_3,\mathbf{v}_1,\bm{\varepsilon})\\
		&\tilde\alpha_1(\mathbf{C},\mathbf{v}_1,\mathbf{v}_2,\mathbf{v}_3,\bm{\varepsilon}) = \tilde\alpha_3(\mathbf{C},\mathbf{v}_2,\mathbf{v}_3,\mathbf{v}_1,\bm{\varepsilon})\\ 
		&\tilde\alpha_2(\mathbf{C},\mathbf{v}_1,\mathbf{v}_2,\mathbf{v}_3,\bm{\varepsilon}) = \tilde\alpha_1(\mathbf{C},\mathbf{v}_2,\mathbf{v}_3,\mathbf{v}_1,\bm{\varepsilon})\\
		&\tilde\alpha_4(\mathbf{C},\mathbf{v}_1,\mathbf{v}_2,\mathbf{v}_3,\bm{\varepsilon}) = \tilde\alpha_4(\mathbf{C},\mathbf{v}_2,\mathbf{v}_3,\mathbf{v}_1,\bm{\varepsilon})\\ 
	\end{aligned}
\end{align}
Note that a redundant constraint $\tilde\alpha_3(\mathbf{C},\mathbf{v}_1,\mathbf{v}_2,\mathbf{v}_3,\bm{\varepsilon}) = \tilde\alpha_2(\mathbf{C},\mathbf{v}_2,\mathbf{v}_3,\mathbf{v}_1,\bm{\varepsilon})$ is removed in (\ref{eq constraint_C3}).

To find the representation for a scalar function $\psi(\mathbf C)$, we can first define an isotropic extension $\hat \psi(\mathbf C, \mathbf v_1, \mathbf v_2, \mathbf v_3, \bm{\varepsilon})$ and then employ (\ref{eq basis_C3}), as
\begin{align}\label{eq psi_C3}
	\begin{array}{l}
		\psi(\mathbf C) = \hat \psi(\mathbf C, \mathbf v_1, \mathbf v_2, \mathbf v_3, \bm{\varepsilon})
		= \hat \psi( \mathbf v_1.\mathbf{C} \mathbf v_1, \mathbf v_2.\mathbf{C} \mathbf v_2, \mathbf v_3.\mathbf{C} \mathbf v_3)
	\end{array}
\end{align}
The additional constraint on (\ref{eq psi_C3}) follows a similar form as (\ref{eq 5.8}) and results in
\begin{align}\label{eq psi_constraint_C3}
	\begin{array}{l}
		\hat \psi( \mathbf v_1.\mathbf{C} \mathbf v_1, \mathbf v_2.\mathbf{C} \mathbf v_2, \mathbf v_3.\mathbf{C} \mathbf v_3)
		= \hat \psi( \mathbf v_2.\mathbf{C} \mathbf v_2, \mathbf v_3.\mathbf{C} \mathbf v_3, \mathbf v_1.\mathbf{C} \mathbf v_1)
	\end{array}
\end{align}

\subsection{Group $\mathcal C_{3v}$}
The symmetry of the point group $\mathcal C_{3v}$ is illustrated in Figure \ref{figure 2}(b). In addition to the 3-fold rotation symmetry, $\mathcal C_{3v}$ also exhibits three reflection planes. The structural tensor set of $\mathcal C_{3v}$ is chosen as $ \{\mathbf{v}_1,\mathbf{v}_2$,$\mathbf{v}_3 \} $ with the structural vectors defined in (\ref{eq structural_C3}). The group $\mathcal C_{3v}$ has two generators $\mathcal{G}^* = \{\mathbf{Q}^{2\pi /3},\bm \sigma_1\}$ (see Table S11 in Supplementary Information). Herein, $\bm \sigma_1 =\begin{bmatrix}
-1 & 0 \\
0 & 1
\end{bmatrix} $ represents the vertical reflection symmetry. Under symmetry transformation $\mathbf{Q}^{2\pi /3}$ and $\bm \sigma_1$, the three structural vectors transform in the following way, as   
\begin{align} \label{eq g*_C3v}
	\begin{array}{ll}
		\mathbf{Q}^{2\pi /3}\mathbf v_1 = \mathbf v_2, \quad &\bm{\sigma}_1\mathbf v_1 = \mathbf v_1 \\
		\mathbf{Q}^{2\pi /3}\mathbf v_2 = \mathbf v_3, \quad & \bm{\sigma}_1\mathbf v_2 = \mathbf v_3 \\
		\mathbf{Q}^{2\pi /3}\mathbf v_3 = \mathbf v_1, \quad & \bm{\sigma}_1\mathbf v_3 =\mathbf v_2 \\
	\end{array}
\end{align}
As it is obvious from (\ref{eq g*_C3v}), $\mathbf{Q}^{2\pi /3}$ permutes all the structural vectors, while ${\bm \sigma _1}$ keeps $\mathbf v_1$ invariant and permutes $\mathbf v_2$ and $\mathbf v_3$. Hence, additional constraints (\ref{eq_constraint}) will be imposed to the representations. 

The representation of tensor functions for $\mathcal{C}_{3v}$ is simpler than $\mathcal C_3$. Firstly, an isotropic extension of $\mathbf{T}(\mathbf C)$ is written as 
\begin{align}\label{eq 5.12}
	\mathbf{T(C)} = \mathbf{\hat T}(\mathbf{C},{\mathbf v_1},{\mathbf v_2},{\mathbf v_3})    
\end{align} 					
By eliminating all $\bm \varepsilon$ terms in (\ref{eq generators_C3}), the representation of (\ref{eq 5.12}) is as follows.
\begin{align} \label{eq generators_redundant_C3v}
	\begin{aligned} & \mathbf{T(C)} = \mathbf{\hat T}(\mathbf{C},{\mathbf v_1},{\mathbf v_2},{\mathbf v_3}) =\alpha_0 \mathbf C +  \alpha_1 (\mathbf v_1 \otimes \mathbf v_1)+ \alpha_2(\mathbf v_2 \otimes \mathbf v_2) + \alpha _3(\mathbf v_3 \otimes \mathbf v_3)
	\end{aligned}
\end{align}
where 
\begin{align}  \label{eq basis_redundant_C3v}
	{\alpha _i = \alpha _i( \mathbf v_1.\mathbf{C}\mathbf v_1, \mathbf v_2.\mathbf{C}\mathbf v_2, \mathbf v_3.\mathbf{C}\mathbf v_3)=\Tilde{\alpha}_i(\mathbf{C}, \mathbf v_1, \mathbf v_2, \mathbf v_3)}, \ i=0,1,2,3  
\end{align} 
Next, additional constraints (\ref{eq_constraint}) must be imposed to (\ref{eq generators_redundant_C3v}). Because this point group has two generators $\mathcal{G}^* = \{\mathbf{Q}^{2\pi /3},\bm \sigma_1\}$, we must impose the constraints for both of them. For $\mathbf{Q} = \mathbf{Q}^{2\pi /3}$, we obtain
\begin{align} \label{eq constraint1_C3v}
	\begin{array}{l}
		\tilde{\alpha} _0(\mathbf{C},\mathbf v_1, \mathbf v_2, \mathbf v_3) = \tilde{\alpha}_0(\mathbf{C},{\mathbf v_2},{\mathbf v_3},{\mathbf v_1})\\
		\tilde{\alpha} _1(\mathbf{C},\mathbf v_1, \mathbf v_2, \mathbf v_3) = \tilde{\alpha}_3(\mathbf{C},{\mathbf v_2},{\mathbf v_3},{\mathbf v_1})\\ 
		\tilde{\alpha}_2(\mathbf{C},\mathbf v_1, \mathbf v_2, \mathbf v_3) = \tilde{\alpha}_1(\mathbf{C},{\mathbf v_2},{\mathbf v_3},{\mathbf v_1})
	\end{array}
\end{align}
and for $\mathbf{Q} = {\bm \sigma _1}$, we obtain
\begin{align} \label{eq constraint2_C3v}
	\begin{array}{ll} 
		\tilde{\alpha}_3(\mathbf{C},\mathbf v_1, \mathbf v_2, \mathbf v_3) = \tilde{\alpha}_2(\mathbf{C},{\mathbf v_1},{\mathbf v_3},{\mathbf v_2})\\ 
	\end{array}
\end{align}
Note that we have removed three redundant constraints in (\ref{eq constraint1_C3v}) and (\ref{eq constraint2_C3v}).

Using the same approach, the representation of a scalar function $\psi(\mathbf C)$ will be as follows.
\begin{align}\label{eq psi_C3v}
	\begin{array}{l}
		\psi(\mathbf C) = \hat \psi(\mathbf C, \mathbf v_1, \mathbf v_2, \mathbf v_3)
		= \hat \psi( \mathbf v_1.\mathbf{C} \mathbf v_1, \mathbf v_2.\mathbf{C} \mathbf v_2, \mathbf v_3.\mathbf{C} \mathbf v_3)
	\end{array}
\end{align}
In addition, the additional constraints for $\mathcal{G}^* = \{\mathbf{Q}^{2\pi /3},\bm \sigma_1\}$ are
\begin{align}\label{eq psi_constraint_C3v}
	\begin{array}{l}
		\hat \psi( \mathbf v_1.\mathbf{C} \mathbf v_1, \mathbf v_2.\mathbf{C} \mathbf v_2, \mathbf v_3.\mathbf{C} \mathbf v_3) \\ = \hat \psi( \mathbf v_2.\mathbf{C} \mathbf v_2, \mathbf v_3.\mathbf{C} \mathbf v_3, \mathbf v_1.\mathbf{C} \mathbf v_1) \\ =\hat \psi( \mathbf v_1.\mathbf{C} \mathbf v_1, \mathbf v_3.\mathbf{C} \mathbf v_3, \mathbf v_2.\mathbf{C} \mathbf v_2)
	\end{array}
\end{align}
By comparing (\ref{eq psi_constraint_C3}) and (\ref{eq psi_constraint_C3v}), we find that the scalar function $\psi (\mathbf C)$ has identical representation for $\mathcal C_3$ and $\mathcal C_{3v}$ but different constraints.

\section{Representation of tensor functions for $\mathcal C_6$ and $\mathcal C_{6v}$} \label{sec:C6}
\subsection{Group $\mathcal C_6$} 
\begin{figure}[!h]
	\centering\includegraphics[width=4 in]{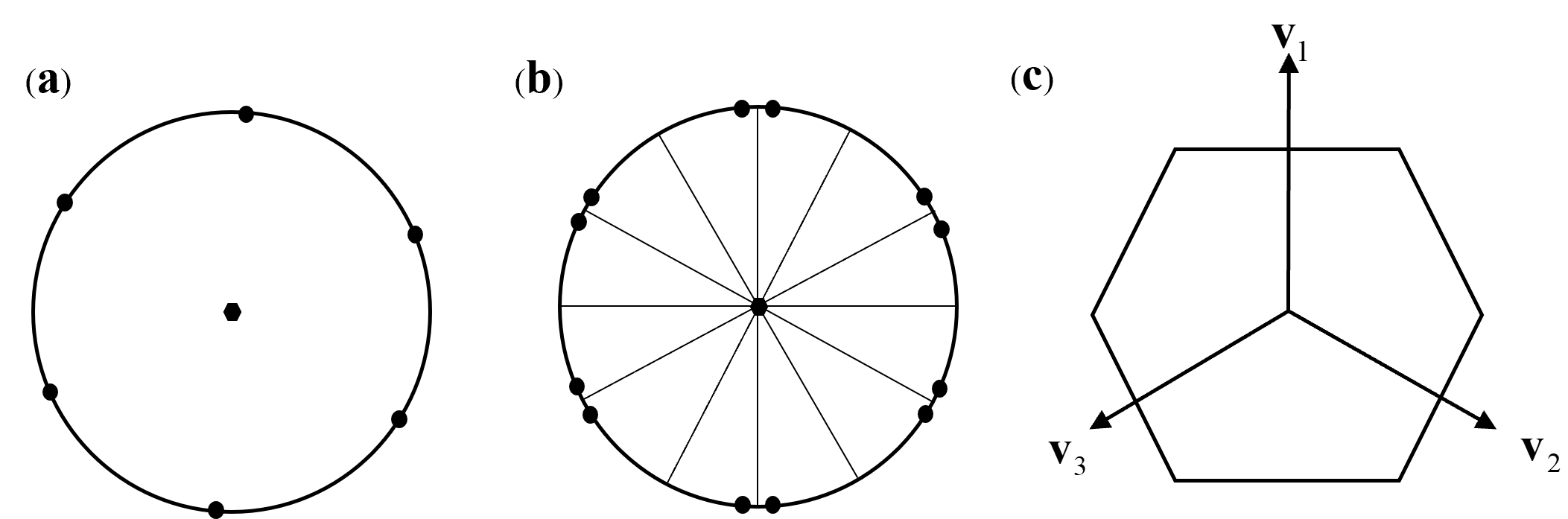}
	\centering \caption{Graphical representation of the 2D point group (a) $\mathcal C_6$, (b) $\mathcal C_{6v}$ and (c) illustration of vectors used to define structural tensor set for point groups $\mathcal C_6$ and $\mathcal C_{6v}$. }
	\label{figure 3}
\end{figure}

As illustrated in Figure \ref{figure 3}(a), the point group $\mathcal C_6$ has 6-fold rotation symmetry and no reflection symmetry. The group $\mathcal C_6$ is a centrosymmetric group, which can be obtained by adding an inversion symmetry to $\mathcal C_3$. The structural tensor set of $\mathcal C_6$ is chosen as $\{\mathbf{M}_1,\mathbf{M}_2, \mathbf{M}_3, \bm \varepsilon \}$, where $\mathbf M_i = \mathbf v_i \otimes \mathbf v_i \ (i=1,2,3)$ with $\mathbf v_i$ defined in (\ref{eq structural_C3}). These four structural tensors have the following matrix form  
\begin{align} \label{eq structural_C6}
	\mathbf{M}_{1}=\begin{bmatrix}0 & 0\\
		0 & 1
	\end{bmatrix},\mathbf{M}_{2}=\left[\begin{array}{cc}
		3/4 & -\sqrt{3}/4\\
		-\sqrt{3}/4 & 1/4
	\end{array}\right],\mathbf{M}_{3}=\left[\begin{array}{cc}
		3/4 & \sqrt{3}/4\\
		\sqrt{3}/4 & 1/4
	\end{array}\right],\bm{\varepsilon}=\begin{bmatrix}0 & 1\\
		-1 & 0
	\end{bmatrix}
\end{align}

The point group $\mathcal C_6$ has one generator $\mathcal{G}^* = \{\mathbf{Q}^{\pi /3}\}$, where $\mathbf{Q}^{\pi /3} =\begin{bmatrix}
1/2 & \sqrt{3}/2 \\
-\sqrt{3}/2 & 1/2
\end{bmatrix} $  represents $\pi/3$-rotation symmetry (see Table S11 in Supplementary Information). Under rotation $\mathbf{Q}^{\pi /3}$, the four structural tensors transform in the following way.   
\begin{align} \label{eq g*_C6}
	\begin{array}{ll}
		\mathbf{Q}^{\pi/3}\mathbf{M}_{1}(\mathbf{Q}^{\pi/3})^{T}=\mathbf{M}_{3}\\
		\mathbf{Q}^{\pi/3}\mathbf{M}_{2}(\mathbf{Q}^{\pi/3})^{T}=\mathbf{M}_{1}\\
		\mathbf{Q}^{\pi/3}\mathbf{M}_{3}(\mathbf{Q}^{\pi/3})^{T}=\mathbf{M}_{2}\\
		\mathbf{Q}^{\pi/3}\bm{\varepsilon}(\mathbf{Q}^{\pi/3})^{T}=\bm{\varepsilon}
	\end{array}
\end{align}
As it is obvious from (\ref{eq g*_C6}), $\mathbf{Q}^{\pi/3}$ permutes
$\mathbf{M}_{1}$, $\mathbf{M}_{2}$, and $\mathbf{M}_{3}$ to each other and keeps $\bm{\varepsilon}$
invariant. Hence, additional constraints must be imposed to the representations.  

The representation of a tensor-valued function $\mathbf{T}(\mathbf C)$ is derived first. Using the structural tensor set, the isotropic extension $\hat{\mathbf T}$ is defined as
\begin{align} \label{eq 6.16}
	\mathbf{T(C)=\hat{T}(C,M_{1},M_{2},M_{3},\bm{\varepsilon})}
\end{align}
Considering that $\mathbf{{M_{1}},{M_{2}},{M_{3}}}$ are symmetric and $\bm{\varepsilon}$ is skew-symmetric $2^{nd}$-order tensors, the representation of (\ref{eq 6.16}) can be obtained using Tables \ref{Table 1} and \ref{Table 2}, as
\begin{align} \label{eq generator_redundant_C6}
	\begin{array}{l}
		\mathbf{T(C)} = \mathbf{\hat T}(\mathbf{C},{\mathbf{M}_1},{\mathbf{M}_2},{\mathbf{M}_3},\bm\varepsilon) = {\alpha _0}\mathbf I + {\alpha _1}\mathbf{C} + {\alpha _2}{\mathbf{M}_1} + {\alpha _3}{\mathbf{M}_2} + {\alpha _4}{\mathbf{M}_3}\\
		\quad \quad \ \ 
		+\alpha_{5}(\mathbf C\bm{\varepsilon}-\bm{\varepsilon}\mathbf{C})+\alpha_{6}(\mathbf M_{1}\bm{\varepsilon}-\bm{\varepsilon}\mathbf M_{1})+\alpha_{7}(\mathbf M_{2}\bm{\varepsilon}-\bm{\varepsilon}\mathbf M_{2})+ \alpha_{8}(\mathbf M_{3}\bm{\varepsilon}-\bm{\varepsilon}\mathbf M_{3})
	\end{array}
\end{align}
where
\begin{align} \label{eq basis_redundant_C6}
	\begin{array}{l}
		{\alpha _i} = {\alpha _i}({{tr}}\mathbf{C},{{ tr}}{\mathbf{C}^2},{{ tr}}{\mathbf{M}_1},{{ tr}}\mathbf{M}_1^2,{{ tr}}{\mathbf{M}_2},{{tr}}\mathbf{M}_2^2,{{tr}}{\mathbf{M}_3},{{ tr}}\mathbf{M}_3^2,{{ tr(}}\mathbf{C}{\mathbf{M}_1}),\\ 
		\quad \quad
		{{tr(}}\mathbf{C}{\mathbf{M}_2}),{{ tr(}}\mathbf{C}{\mathbf{M}_3}),{{tr(}}{\mathbf{M}_1}{\mathbf{M}_2}),{{tr(}}{\mathbf{M}_1}{\mathbf{M}_3}),{{ tr(}}{\mathbf{M}_2}{\mathbf{M}_3}), tr\bm{\varepsilon}^2,\\
		\quad \quad 
		tr(\mathbf C \mathbf M_1 \bm{\varepsilon}), tr(\mathbf C \mathbf M_2 \bm{\varepsilon}), 	tr(\mathbf C \mathbf M_3 \bm{\varepsilon}),	tr(\mathbf M_1 \mathbf M_2 \bm{\varepsilon}), tr(\mathbf M_1 \mathbf M_3 \bm{\varepsilon}),\\
		\quad \quad
		tr(\mathbf M_2 \mathbf M_3 \bm{\varepsilon}))
	\end{array}
\end{align}
After eliminating the redundant terms in (\ref{eq generator_redundant_C6}) and (\ref{eq basis_redundant_C6}), the representation of $\mathbf T(\mathbf C)$ is obtained as
\begin{align} \label{eq  generator_C6}
	\begin{array}{ll}
		\mathbf T(\mathbf C)=\hat{\mathbf T}(\mathbf C,\mathbf M_{1},\mathbf M_{2},\mathbf M_{3},\bm{\varepsilon})=\alpha_0 \mathbf C + \alpha_{1}\mathbf M_{1}+\alpha_{2}\mathbf M_{2}+\alpha_{3}\mathbf M_{3} +\alpha_{4}(\mathbf C\bm{\varepsilon}-\bm{\varepsilon}\mathbf{C})
	\end{array}
\end{align}
and 
\begin{align} \label{eq  basis_C6}
	\begin{array}{ll}
		\alpha_{i}=\alpha_{i}(tr(\mathbf{C}\mathbf{M}_{1}),tr(\mathbf{C}\mathbf{M}_{2}),tr(\mathbf{C}\mathbf{M}_{3}))\\
		\quad \ =\tilde{\alpha}_i(\mathbf{C}, \mathbf M_1, \mathbf M_2, \mathbf M_3, \bm{\varepsilon}); \quad i=0,1,2,3,4
	\end{array}
\end{align}
Additional constraints (\ref{eq_constraint}) must be imposed to (\ref{eq  generator_C6}) for $\mathcal{G}^* = \{\mathbf{Q}^{\pi /3}\}$. Considering $\mathbf{Q} = \mathbf{Q}^{\pi /3}$, the scalar coefficient functions $\tilde{\alpha}_i$ must satisfy the following constraints. 
\begin{align} \label{eq constraint_C6}
	\begin{array} {ll}
		\tilde\alpha_0(\mathbf{C},\mathbf{M}_1,\mathbf{M}_2,\mathbf{M}_3,\bm{\varepsilon}) = \tilde\alpha_0(\mathbf{C},\mathbf{M}_3,\mathbf M_1,\mathbf M_2,\bm{\varepsilon})\\
		\tilde\alpha_1(\mathbf{C},\mathbf{M}_1,\mathbf{M}_2,\mathbf{M}_3,\bm{\varepsilon}) = \tilde\alpha_2(\mathbf{C},\mathbf{M}_3,\mathbf M_1,\mathbf M_2,\bm{\varepsilon})\\ 
		\tilde\alpha_2(\mathbf{C},\mathbf{M}_1,\mathbf{M}_2,\mathbf{M}_3,\bm{\varepsilon}) = \tilde\alpha_3(\mathbf{C},\mathbf{M}_3,\mathbf M_1,\mathbf M_2,\bm{\varepsilon})\\ 
		\tilde\alpha_4(\mathbf{C},\mathbf{M}_1,\mathbf{M}_2,\mathbf{M}_3,\bm{\varepsilon}) = \tilde\alpha_4(\mathbf{C},\mathbf{M}_3,\mathbf M_1,\mathbf M_2,\bm{\varepsilon})\\ 
	\end{array}
\end{align}
Note that a redundant constraint $\tilde\alpha_3(\mathbf{C},\mathbf{M}_1,\mathbf{M}_2,\mathbf{M}_3,\bm{\varepsilon}) = \tilde\alpha_1(\mathbf{C},\mathbf{M}_3,\mathbf M_1,\mathbf M_2,\bm{\varepsilon})$ is removed in (\ref{eq constraint_C6}).

The representation of a scalar function $\psi(\mathbf C)$ follows a similar form of (\ref{eq basis_C6}) as
\begin{align}\label{eq psi_C6}
	\begin{array}{l}
		\psi(\mathbf C) = \hat \psi(\mathbf C, \mathbf M_1, \mathbf M_2, \mathbf M_3, \bm{\varepsilon})
		= \hat \psi(tr(\mathbf{C}\mathbf{M}_{1}),tr(\mathbf{C}\mathbf{M}_{2}),tr(\mathbf{C}\mathbf{M}_{3}))
	\end{array}
\end{align}
In addition, the additional constraint is
\begin{align}\label{eq psi_constraint_C6}
	\begin{array}{l}
		\hat \psi(tr(\mathbf{C}\mathbf{M}_{1}),tr(\mathbf{C}\mathbf{M}_{2}),tr(\mathbf{C}\mathbf{M}_{3}))
		= \hat \psi(tr(\mathbf{C}\mathbf{M}_{3}),tr(\mathbf{C}\mathbf{M}_{1}),tr(\mathbf{C}\mathbf{M}_{2}))
	\end{array}
\end{align}

\subsection{Group $\mathcal C_{6v}$} 
As illustrated in Figure \ref{figure 3}(b), the point group $\mathcal C_{6v}$ has 6-fold rotation symmetry and reflection symmetry. The structural tensor set of $\mathcal C_{6v}$ is chosen as $\{\mathbf M_1, \mathbf M_2, \mathbf M_3\}$ with $\mathbf M_i \ (i=1,2,3)$ defined in (\ref{eq structural_C6}). The group generators of $\mathcal C_{6v}$ are $\mathcal{G}^* = \{\mathbf{Q}^{\pi /3},\bm \sigma_{d_1}\}$, where $\mathbf{Q}^{\pi /3}$ is $\pi/3$-rotation and $\bm \sigma_{d_1} = \begin{bmatrix}
-1 & 0 \\
0 & 1
\end{bmatrix}$ is vertical reflection (see Table S11 in Supplementary Information). Under $\mathbf{Q}^{\pi /3}$ and $\bm \sigma_{d_1}$, the three structural tensors transform in the following way.   
\begin{align} \label{eq g*_C6v}
	\begin{array}{l}
		\mathbf{Q}^{\pi /3}{\mathbf{M}_1}{(\mathbf{Q}^{\pi /3})^T} = {\mathbf{M}_3},\\
		\mathbf{Q}^{\pi /3}{\mathbf{M}_2}{(\mathbf{Q}^{\pi /3})^T} = {\mathbf{M}_1},\\
		\mathbf{Q}^{\pi /3}{\mathbf{M}_3}{(\mathbf{Q}^{\pi /3})^T} = {\mathbf{M}_2},
	\end{array}  \begin{array}{l}
		{\bm \sigma _{{d_1}}}{\mathbf{M}_1}{\bm \sigma _{{d_1}}^T} = {\mathbf{M}_1}\\
		{\bm \sigma _{{d_1}}}{\mathbf{M}_2}{\bm \sigma _{{d_1}}^T} = {\mathbf{M}_3}\\
		{\bm \sigma _{{d_1}}}{\mathbf{M}_3}{\bm \sigma _{{d_1}}^T} = {\mathbf{M}_2}
	\end{array}	    
\end{align}
As it is obvious from (\ref{eq  g*_C6v}), $\mathbf{Q}^{\pi /3}$ permutes all the structural tensors, $\bm \sigma_{d_1}$ keeps $\mathbf{M}_1$ invariant and permutes $\mathbf{M}_2$ to $\mathbf{M}_3$. Thus, additional constraints must be imposed to the representations.

The representation of a tensor valued function $\mathbf{T(C)}$ is presented first. Using the structural tensor set, we can define an isotropic extension $\hat{\mathbf T}$ as  
\begin{align} \label{eq 6.4}
	\mathbf{T(C)} = \mathbf{\hat T}(\mathbf{C},{\mathbf{M}_1},{\mathbf{M}_2},{\mathbf{M}_3})    
\end{align} 					
Considering that all the arguments are symmetric, the representation of (\ref{eq 6.4}) can be obtained using Tables \ref{Table 1} and \ref{Table 2}, as
\begin{align} \label{eq generator_redundant_C6v}
	\mathbf{T(C)} = \mathbf{\hat T}(\mathbf{C},{\mathbf{M}_1},{\mathbf{M}_2},{\mathbf{M}_3}) = {\alpha _0}\mathbf I + {\alpha _1}\mathbf{C} + {\alpha _2}{\mathbf{M}_1} + {\alpha _3}{\mathbf{M}_2} + {\alpha _4}{\mathbf{M}_3}
\end{align}
where 
\begin{align} \label{eq basis_redundant_C6v}
	\begin{array}{l}
		{\alpha _i} = {\alpha _i}({{tr}}\mathbf{C},{{ tr}}{\mathbf{C}^2},{{ tr}}{\mathbf{M}_1},{{ tr}}\mathbf{M}_1^2,{{ tr}}{\mathbf{M}_2},{{tr}}\mathbf{M}_2^2,{{tr}}{\mathbf{M}_3},{{ tr}}\mathbf{M}_3^2,{{ tr(}}\mathbf{C}{\mathbf{M}_1}),\\ \quad 
		{{tr(}}\mathbf{C}{\mathbf{M}_2}),{{ tr(}}\mathbf{C}{\mathbf{M}_3}),{{tr(}}{\mathbf{M}_1}{\mathbf{M}_2}),{{tr(}}{\mathbf{M}_1}{\mathbf{M}_3}),{{ tr(}}{\mathbf{M}_2}{\mathbf{M}_3}))
	\end{array}
\end{align}
Then, after eliminating redundant terms in (\ref{eq generator_redundant_C6v}) and (\ref{eq basis_redundant_C6v}), the representation of $\mathbf T (\mathbf C )$ is finally expressed as follows. 
\begin{align} \label{eq generator_C6v}
	\mathbf{T(C)} = \mathbf{\hat T}(\mathbf{C},{\mathbf{M}_1},{\mathbf{M}_2},{\mathbf{M}_3}) =\alpha_0 \mathbf C +   {\alpha _1}{\mathbf{M}_1} + {\alpha _2}{\mathbf{M}_2} + {\alpha _3}{\mathbf{M}_3}
\end{align}
where 
\begin{align}  \label{eq basis_C6v}
	\begin{array}{l}
		{\alpha _i} = {\alpha _i}({{ tr(}}\mathbf{C}{\mathbf{M}_1}),{{tr(}}\mathbf{C}{\mathbf{M}_2}),{{ tr(}}\mathbf{C}{\mathbf{M}_3}))  \\ \quad \
		=\tilde{\alpha}_i(\mathbf{C},{\mathbf{M}_1},{\mathbf{M}_2},{\mathbf{M}_3}) ; \quad i=0,1,2,3 
	\end{array}
\end{align}
Additional constraints (\ref{eq_constraint}) must be imposed to the representation (\ref{eq generator_C6v}) for group generators $\mathcal{G}^* = \{\mathbf{Q}^{\pi /3},\bm \sigma_{d_1}\}$. Specifically, for $\mathbf{Q}^{\pi /3}$, this requires $\mathbf{\hat T(C},{\mathbf{M}_1},{\mathbf{M}_2},{\mathbf{M}_3}) = \mathbf{\mathbf{\hat T}(\mathbf{C}},{\mathbf{M}_3},{\mathbf{M}_1},{\mathbf{M}_2})$; and for $\bm \sigma_{d_1}$, this requires $\mathbf{\hat T(C},{\mathbf{M}_1},{\mathbf{M}_2},{\mathbf{M}_3}) = \mathbf{\mathbf{\hat T}(\mathbf{C}},{\mathbf{M}_1},{\mathbf{M}_3},{\mathbf{M}_2})$. Further, we can find the constraints to the scalar coefficient functions. 
For $\mathbf{Q} = \mathbf{Q}^{\pi /3}$, we obtain
\begin{align} \label{eq constraint1_C6v}
	\begin{array}{ll}
		\tilde\alpha _0(\mathbf{C},{\mathbf{M}_1},{\mathbf{M}_2},{\mathbf{M}_3}) = \tilde\alpha _0(\mathbf{C},{\mathbf{M}_3},{\mathbf{M}_1},{\mathbf{M}_2}){\rm{ }}\\ 
		\tilde\alpha _1(\mathbf{C},{\mathbf{M}_1},{\mathbf{M}_2},{\mathbf{M}_3}) = \tilde\alpha _2(\mathbf{C},{\mathbf{M}_3},{\mathbf{M}_1},{\mathbf{M}_2}){\rm{ }}\\ 
		\tilde\alpha _3(\mathbf{C},{\mathbf{M}_1},{\mathbf{M}_2},{\mathbf{M}_3}) = \tilde\alpha _1(\mathbf{C},{\mathbf{M}_3},{\mathbf{M}_1},{\mathbf{M}_2}){\rm{ }}\\ 
	\end{array}
\end{align}
and for $\mathbf{Q} = {\bm \sigma _{{d_1}}}$, we obtain
\begin{align} \label{eqconstraint2_C6v}
	\begin{array}{ll}
		\tilde\alpha _2(\mathbf{C},{\mathbf{M}_1},{\mathbf{M}_2},{\mathbf{M}_3}) = \tilde\alpha _3(\mathbf{C},{\mathbf{M}_1},{\mathbf{M}_3},{\mathbf{M}_2})
	\end{array}
\end{align}
Note that we have removed three redundant constraints in (\ref{eq constraint1_C6v}) and (\ref{eqconstraint2_C6v}).

The representation of a scalar function $\psi(\mathbf C)$ is similar to (\ref{eq basis_C6v}), as
\begin{align}\label{eq psi_C6v}
	\begin{array}{l}
		\psi(\mathbf C) = \hat \psi(\mathbf C, \mathbf M_1, \mathbf M_2, \mathbf M_3)
		= \hat \psi(tr(\mathbf{C}\mathbf{M}_{1}),tr(\mathbf{C}\mathbf{M}_{2}),tr(\mathbf{C}\mathbf{M}_{3}))
	\end{array}
\end{align}
and the additional constraints are
\begin{align}\label{eq psi_constraint_C6v}
	\begin{array}{l}
		\hat \psi(tr(\mathbf{C}\mathbf{M}_{1}),tr(\mathbf{C}\mathbf{M}_{2}),tr(\mathbf{C}\mathbf{M}_{3}))\\
		= \hat \psi(tr(\mathbf{C}\mathbf{M}_{3}),tr(\mathbf{C}\mathbf{M}_{1}),tr(\mathbf{C}\mathbf{M}_{2}))\\
		=	\hat \psi(tr(\mathbf{C}\mathbf{M}_{1}),tr(\mathbf{C}\mathbf{M}_{3}),tr(\mathbf{C}\mathbf{M}_{2}))
	\end{array}
\end{align}

\section{Further remarks} \label{sec:remarks}

Some remarks are added below regarding the presented theory in general.

\textbf{Remark 7.1.} As mentioned in Section \ref{sec:struct_tensor_set}, each point group has various options for the structural tensor set, and Table \ref{Table 4} only suggests one option. The representation of tensor functions is non-unique and dependent on the chosen structural tensor set. 

\textbf{Remark 7.2.} For the representation of scalar functions or second-order symmetric tensor functions, the results for the point groups $\mathcal C_{6v}$ and $\mathcal C_{3v}$ are equivalent. This comes from the fact that these functions exhibit a centrosymmetric feature so their symmetry includes both structural symmetry from the point group and the inherent inversion symmetry of the functions. Thus in this case, the results for several 2D point group pairs are equivalent: $\mathcal C_{6v} \sim \mathcal C_{3v}$, $\mathcal C_{6} \sim \mathcal C_{3}$, $\mathcal C_{2} \sim \mathcal C_{1}$ and $\mathcal C_{2v} \sim \mathcal C_{1v}$. In short, for centrosymmetric properties, we only need to study the corresponding Laue groups.

\textbf{Remark 7.3.} In this work, we only present the representations of scalar- and symmetric tensor-valued functions for each point group. For skew-symmetric tensor-valued functions, one may follow a similar procedure to derive the representations. They are not used very often so we will not study them in this work. Moreover, it is worth mentioning that, using the approach presented in this work, the elastic strain energy can be formulated as a scalar-valued tensor function. The elasticity tensor can then be derived by taking the second-order derivative of the strain energy with respect to the strain tensor. In addition, an important related approach not explored here is the spectral decomposition of the elasticity tensor, as utilized in linear elasticity. This idea, originally introduced by Kelvin and recently revisited in works such as \cite{kowalczyk2009review}, allows the representation of the elasticity tensor through its spectral form and can be naturally extended to general tensor functions, although the resulting decomposition may involve asymmetric components. This direction provides an alternative and insightful perspective that may complement the present formulation. 

\textbf{Remark 7.4.} In this work, we focus on the representation of tensor functions for solids and their associated symmetry groups. While we recognize the importance of other material symmetries, such as the full unimodular group relevant to fluids and its subgroups for fluid crystals. However, addressing these symmetries falls outside the scope of this work. It's worth noting that by anisotropic we mean not-isotropic in $\mathcal O(2)$.

\textbf{Remark 7.5.} It is worth noting that tensor-valued functions such as $\mathbf{T}(\mathbf{v}, \mathbf{A}, \mathbf{W})$ can also be derived from scalar-valued functions $\Psi(\mathbf{v}, \mathbf{A}, \mathbf{W}, \mathbf{x}, \mathbf{X})=\Psi(I_k)$, where $I_k$ denotes scalar invariants (integrity basis) constructed from the arguments $\{\mathbf{v}, \mathbf{A}, \mathbf{W}, \mathbf{x}, \mathbf{X}\}$, with $\mathbf{x}$ and $\mathbf{X}$ denoting auxiliary vector and tensor arguments, respectively. By applying the chain rule, one obtains a representation of the form
\[
\mathbf{T}(\mathbf{v}, \mathbf{A}, \mathbf{W}) = \sum_j \alpha_j \mathbf{G}_j, \quad \text{where } \alpha_j = \frac{\partial \Psi}{\partial I_j}, \quad \mathbf{G}_j = \frac{\partial I_j}{\partial \mathbf{X}}.
\]
This approach offers an alternative way to derive tensor-valued functions from scalar functions, although the equations are lengthy due to the additional arguments introduced. The foundational idea of generating tensor functions via derivatives of scalar functions traces back to Rivlin \cite{rivlin_material_1980}. While our current work emphasizes direct construction via structural tensors and symmetry constraints, this alternative approach is certainly complementary and worthy of exploration in future studies.

\begin{figure}[!h]
	\centering\includegraphics[width=5.25in]{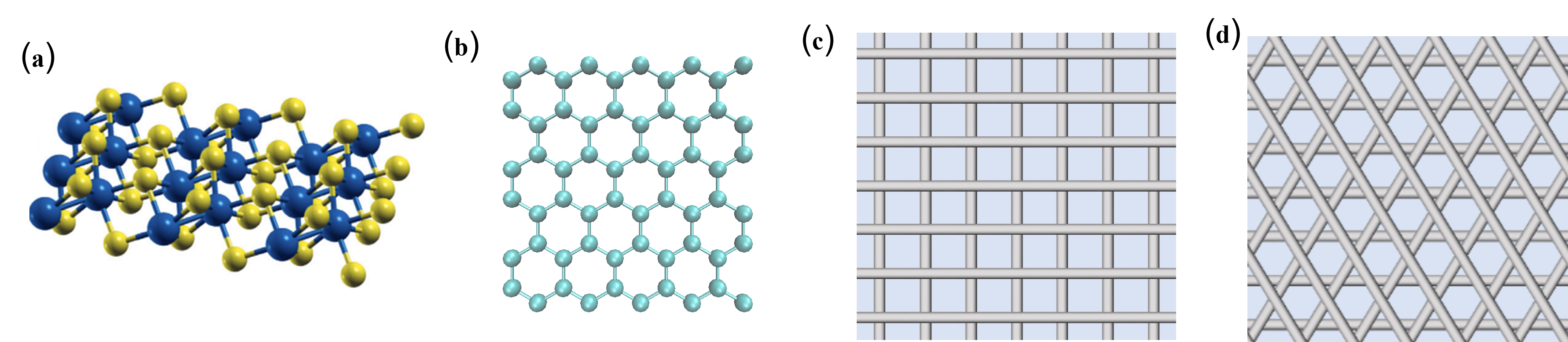}
	\centering \caption{(a) Triclinic 2D rhenium disulfide \cite{maalouf_anisotropic_2023}. (b) Graphene. Fibrous Composites with (c) square and (d) triangular lattice reinforcement.} 
	\label{figure 4}
\end{figure}

\textbf{Remark on applications} Representation of 2D anisotropic scalar- and tensor-valued functions has broad applications in engineering and materials science. 2D materials and structures include 2D nanomaterials, thin films, plates, composite films and laminates, 2D lattice materials, etc. Some typical materials are illustrated in Figure \ref{figure 4} as examples. For scalar-valued functions, the theory present in this work can be used to model hyperelastic strain energy functions of elastomers, soft composites, and biological tissues, as well as yield and failure criteria \cite{boehler_yielding_1977} of 2D materials and structures. On the other hand, the representation of tensor-valued functions enables researchers to model mechanical, physical, and mechanophysical properties of 2D materials \cite{madadi_finite_2024, Rahman2025}, including stress-strain relationship, dielectric property, and electrical/thermal conductivity tensor, etc. Certainly, this work only shows the general forms of the scalar- and tensor-valued functions. The detailed functions need to be devised and fitted using data.

\section{Conclusions} \label{sec:conclusion}
In this study, we present a comprehensive study on the representation theory of tensor functions for 2D point groups using lower-order structural tensors. The representation theory is foundational to constitutive modeling of anisotropic materials. Existing challenge is that some point groups involve higher-order structural tensors and are difficult to use for applications. We adopt a reformulated theory proposed by Man and Goddard to overcome this obstacle. By proposing appropriate "structural tensor set" for all 2D point groups, we derived the representations of scalar- and tensor-valued functions using lower-order structural tensors. Correspondingly, the functional bases, tensor generators, and additional constraints of each point group are provided. We argue that for the six point groups ($\mathcal{C}_1$, $\mathcal{C}_{1v}$, $\mathcal{C}_2$, $\mathcal{C}_{2v}$, $\mathcal{C}_\infty$, $\mathcal{C}_{\infty v}$), the original formulation of Boehler and Liu together with Zheng's structural tensors should be used. In contrast, for the six point groups ($\mathcal C_3$, $\mathcal C_{3v}$, $\mathcal C_4$, $\mathcal C_{4v}$, $\mathcal C_6$, $\mathcal C_{6v}$), the reformulation of Man and Goddard together with our proposed structural tensor sets should be used.  

The representation theory present in this work open new opportunities and practical solutions to the future research of mechanical and physical behaviors of anisotropic materials, including constitutive modeling, failure analysis, electrical and thermal properties, among others. In the future, similar theory should be established for 3D point groups that would have greater impact to the modeling of anisotropic materials. In addition, future works can also be directed towards discovering specific functions for the constitutive modeling by using novel methods such as symbolic artificial intelligence. Recent advances in symbolic artificial intelligence enable the discovery of analytical expressions from data directly. The theory established in this work can provide general tensor formulations and constraints for the symbolic artificial intelligence algorithms. 

\section*{Acknowledgements}
M.M and P.Z. acknowledge the funding support from National Science Foundation (CMMI-2244952, CMMI-2143297). L.C. acknowledges the funding support from National Science Foundation (CMMI-2244953).
\section*{Declaration of Interest Statement}
The authors declare that they have no known competing financial interests or personal relationships that could have appeared to influence the work reported in this paper.

\bibliographystyle{ieeetr}
\bibliography{Anisotropic_2D.bib}

\begin{thebibliography}{10}

\bibitem{ottosen2005mechanics}
N.~S. Ottosen and M.~Ristinmaa, {\em The mechanics of constitutive modeling}.
\newblock Amsterdam, Netherlands: Elsevier, 2005.

\bibitem{zheng_theory_1994}
Q.-S. Zheng, ``Theory of {Representations} for {Tensor} {Functions}—{A}
  {Unified} {Invariant} {Approach} to {Constitutive} {Equations},'' {\em
  Applied Mechanics Reviews}, vol.~47, pp.~545--587, Nov. 1994.

\bibitem{rivlin_large_1948}
R.~S. Rivlin, ``Large elastic deformations of isotropic materials iv. further
  developments of the general theory,'' {\em Philosophical Transactions of the
  Royal Society of London. Series A, Mathematical and Physical Sciences},
  vol.~241, no.~835, pp.~379--397, 1948.

\bibitem{rivlin_stress_1997}
R.~S. Rivlin and J.~L. Ericksen, ``Stress-deformation relations for isotropic
  materials,'' {\em Collected Papers of RS Rivlin: Volume I and II},
  pp.~911--1013, 1997.

\bibitem{rivlin_further_1955}
R.~S. Rivlin, ``Further remarks on the stress-deformation relations for
  isotropic materials,'' in {\em Collected Papers of RS Rivlin: Volume I and
  II}, pp.~1014--1035, Berlin, Germany: Springer, 1955.

\bibitem{pipkin_formulation_1959}
A.~C. Pipkin and R.~Rivlin, ``The formulation of constitutive equations in
  continuum physics. i,'' in {\em Collected Papers of RS Rivlin: Volume I and
  II}, pp.~1111--1126, Berlin, Germany: Springer, 1959.

\bibitem{noll_representations_1970}
W.~Noll, ``Representations of certain isotropic tensor functions,'' {\em Archiv
  der Mathematik}, vol.~21, no.~1, pp.~87--90, 1970.

\bibitem{wang_representations_1969}
C.~Wang, ``On representations for isotropic functions: Part i. isotropic
  functions of symmetric tensors and vectors,'' {\em Archive for Rational
  Mechanics and Analysis}, vol.~33, pp.~249--267, Jan. 1969.

\bibitem{wangb_representations_1969}
C.~Wang, ``On representations for isotropic functions: Part ii. isotropic
  functions of skew-symmetric tensors, symmetric tensors, and vectors,'' {\em
  Archive for Rational Mechanics and Analysis}, vol.~33, pp.~268--287, 1969.

\bibitem{wang_new_1970}
C.~C. Wang, ``A new representation theorem for isotropic functions: An answer
  to professor {G}.{F}. smith's criticism of my papers on representations for
  isotropic functions: Part 1. scalar-valued isotropic functions,'' {\em
  Archive for Rational Mechanics and Analysis}, vol.~36, pp.~166--197, 1970.

\bibitem{wang_corrigendum_1971}
C.~C. Wang, ``Corrigendum to my recent papers on “{Representations} for
  isotropic functions”,'' {\em Archive for Rational Mechanics and Analysis},
  vol.~43, pp.~392--395, 1971.

\bibitem{smith_fundamental_1970}
G.~Smith, ``On a fundamental error in two papers of {C}.-{C}. {Wang} “{On}
  representations for isotropic functions, parts {I} and {II}”,'' {\em
  Archive for Rational Mechanics and Analysis}, vol.~36, pp.~161--165, 1970.

\bibitem{smith_isotropic_1971}
G.~Smith, ``On isotropic functions of symmetric tensors, skew-symmetric tensors
  and vectors,'' {\em International Journal of Engineering Science}, vol.~9,
  no.~10, pp.~899--916, 1971.

\bibitem{boehler_irreducible_1977}
J.~P. Boehler, ``On {irreducible} {representations} for {isotropic} {scalar}
  {functions},'' {\em ZAMM - Journal of Applied Mathematics and Mechanics /
  Zeitschrift für Angewandte Mathematik und Mechanik}, vol.~57, no.~6,
  pp.~323--327, 1977.

\bibitem{boehler_lois_1978}
J.-P. Boehler, ``Lois de comportement anisotrope des milieux continus.,'' {\em
  J. M\'ecanique}, vol.~17, no.~2, pp.~153--190, 1978.

\bibitem{boehler_simple_1979}
J.-P. Boehler, ``A simple derivation of representations for non-polynomial
  constitutive equations in some cases of anisotropy,'' {\em ZAMM-Journal of
  Applied Mathematics and Mechanics/Zeitschrift f{\"u}r Angewandte Mathematik
  und Mechanik}, vol.~59, no.~4, pp.~157--167, 1979.

\bibitem{Boehler1987}
J.-P. Boehler, {\em Applications of tensor functions in solid mechanics}.
\newblock Vienna: Springer, 1987.

\bibitem{spencer_formulation_1982}
A.~Spencer, ``The formulation of constitutive equation for anisotropic
  solids,'' in {\em Mechanical Behavior of Anisotropic Solids/Comportment
  M{\'e}chanique des Solides Anisotropes: Proceedings of the Euromech
  Colloquium 115 Villard-de-Lans, June 19--22, 1979/Colloque Euromech 115
  Villard-de-Lans, 19--22 juin 1979}, pp.~3--26, Springer, Dordrecht, 1982.

\bibitem{liu_representations_1982}
I.-S. Liu {\em et~al.}, ``On representations of anisotropic invariants,'' {\em
  International Journal of Engineering Science}, vol.~20, no.~10,
  pp.~1099--1109, 1982.

\bibitem{betten_irreducible_1987}
J.~Betten, ``Irreducible invariants of fourth-order tensors,'' {\em
  Mathematical Modelling}, vol.~8, pp.~29--33, 1987.

\bibitem{zheng_two-dimensional_1997}
Q.-S. Zheng, ``Two-dimensional tensor function representation for all kinds of
  material symmetry,'' {\em Proceedings of the Royal Society of London. Series
  A: Mathematical and Physical Sciences}, vol.~443, pp.~127--138, Jan. 1993.

\bibitem{zheng_note_1994}
Q.-S. Zheng, ``A {note} on {representation} for {isotropic} {functions} of
  4th-{order} {tensors} in 2-{dimensional} {space},'' {\em ZAMM - Journal of
  Applied Mathematics and Mechanics / Zeitschrift für Angewandte Mathematik
  und Mechanik}, vol.~74, no.~8, pp.~357--359, 1994.

\bibitem{xiao_isotropic_1996}
H.~Xiao, ``On {isotropic} {extension} of {anisotropic} {tensor} {functions},''
  {\em ZAMM - Journal of Applied Mathematics and Mechanics / Zeitschrift für
  Angewandte Mathematik und Mechanik}, vol.~76, no.~4, pp.~205--214, 1996.

\bibitem{man2018remarks}
C.-S. Man and J.~D. Goddard, ``Remarks on isotropic extension of anisotropic
  constitutive functions via structural tensors,'' {\em Mathematics and
  Mechanics of Solids}, vol.~23, no.~4, pp.~554--563, 2018.

\bibitem{de2012structure}
M.~De~Graef and M.~E. McHenry, {\em Structure of Materials: An Introduction to
  Crystallography, Diffraction and Symmetry}.
\newblock Cambridge: Cambridge University Press, 2~ed., 2012.

\bibitem{Itskov2007}
M.~Itskov, {\em Tensor algebra and tensor analysis for engineers}.
\newblock Berlin, Germany: Springer, 2007.

\bibitem{korsgaard_representation_1990}
J.~Korsgaard, ``On the representation of two-dimensional isotropic functions,''
  {\em International Journal of Engineering Science}, vol.~28, pp.~653--662,
  Jan. 1990.

\bibitem{adkins_symmetry_1959}
J.~Adkins, ``Symmetry relations for orthotropic and transversely isotropic
  materials,'' {\em Archive for Rational Mechanics and Analysis}, vol.~4,
  pp.~193--213, 1959.

\bibitem{adkins_further_1960}
J.~Adkins, ``Further symmetry relations for transversely isotropic materials,''
  {\em Archive for Rational Mechanics and Analysis}, vol.~5, pp.~263--274,
  1960.

\bibitem{ajm_theory_1971}
A.~Spencer, {\em Theory of Invariants}.
\newblock New York: Academic Press, 1971.

\bibitem{lokhin_nonlinear_1963}
V.~Lokhin and L.~Sedov, ``Nonlinear tensor functions of several tensor
  arguments,'' {\em Journal of Applied Mathematics and Mechanics}, vol.~27,
  no.~3, pp.~597--629, 1963.

\bibitem{smith_anisotropic_1957}
G.~Smith and R.~S. Rivlin, ``The anisotropic tensors,'' {\em Quarterly of
  Applied Mathematics}, vol.~15, no.~3, pp.~308--314, 1957.

\bibitem{kowalczyk2009review}
K.~Kowalczyk-Gajewska and J.~Ostrowska-Maciejewska, ``Review on spectral
  decomposition of {Hooke}’s tensor for all symmetry groups of linear elastic
  material,'' {\em Engineering Transactions}, vol.~57, no.~3-4, pp.~145--183,
  2009.

\bibitem{rivlin_material_1980}
R.~S. Rivlin, ``Material symmetry and constitutive equations,'' {\em
  Ingenieur-Archiv}, vol.~49, pp.~325--336, Aug. 1980.

\bibitem{boehler_yielding_1977}
J.~Boehler and A.~Sawczuk, ``On yielding of oriented solids,'' {\em Acta
  Mechanica}, vol.~27, no.~1, pp.~185--204, 1977.

\bibitem{madadi_finite_2024}
M.~Madadi and P.~Zhang, ``Finite-size effect on the percolation and
  electromechanical behaviors of liquid metal particulate composites,'' {\em
  Soft Matter}, vol.~20, no.~5, pp.~1061--1069, 2024.

\bibitem{Rahman2025}
A.~Rahman, M.~Madadi, J.~Ma, and P.~Zhang, ``Overcoming
  conductivity-stretchability tradeoff in soft conductive composites through
  liquid metal junctions,'' July 2025.

\bibitem{maalouf_anisotropic_2023}
S.~R. Maalouf and S.~S. Vel, ``Anisotropic elastic properties of triclinic {2D}
  materials using density functional theory with application to rhenium
  disulfide,'' {\em Computational Condensed Matter}, vol.~34, p.~e00790, 2023.

\end{thebibliography}

\newpage


\section*{Supplementary material}

Supplementary information for

# Representation of tensor functions using low-order structural tensors: two-dimensional point groups

Mohammad Madadi[1], Lin Cheng[2] and Pu Zhang[1,*]

[1]Department of Mechanical Engineering, State University of New York at Binghamton, Binghamton, NY 13902, United States

[2]Department of Mechanical Engineering, University of Maryland, College Park, MD 20742, United States

**List of Supplementary Documents:**

Figure S1. Symmetry hierarchy of two-dimensional crystallographic point groups.
Figure S2. Graphical representation of the ten 2D crystallographic point groups.
Figure S3. The representation of mirror symmetry operators for point group $\mathcal{C}_{4v}$
Figure S4. The representation of mirror symmetry operators for point group $\mathcal{C}_{3v}$.
Figure S5. The representation of mirror symmetry operators for point group $\mathcal{C}_{6v}$.

Table S1. Multiplication table of the point group $\mathcal{C}_1$.
Table S2. Multiplication table of the point group $\mathcal{C}_2$.
Table S3. Multiplication table of the point group $\mathcal{C}_{1v}$.
Table S4. Multiplication table of the point group $\mathcal{C}_{2v}$.
Table S5. Multiplication table of the point group $\mathcal{C}_4$.
Table S6. Multiplication table of the point group $\mathcal{C}_{4v}$.
Table S7. Multiplication table of the point group $\mathcal{C}_3$.
Table S8. Multiplication table of the point group $\mathcal{C}_{3v}$.
Table S9. Multiplication table of the point group $\mathcal{C}_6$.
Table S10. Multiplication table of the point group $\mathcal{C}_{6v}$.
Table S11. Members and generators for the two-dimensional crystallographic and compact point groups.

[*]Corresponding author. Email: pzhang@binghamton.edu

## Multiplication table for 2D point groups:

The multiplication tables for all crystallographic point groups in two-dimensional space are provided in Table S1. to Table S10. In the Tables, $\mathcal{G}$ is the set of all members of the point group. In the tables, in a two-dimensional Euclidean space $\mathbf{R}^2$, rotations through $\alpha = 2\pi/n$, $n \in N$, around the axis a with $\|\mathbf{a}\|=1$ and normal to the lattice plane are defined by

$$\mathbf{Q}^\alpha = \begin{pmatrix} \cos\alpha & \sin\alpha \\ -\sin\alpha & \cos\alpha \end{pmatrix} \qquad (1)$$

In addition, symmetry hierarchy and graphical representation of the 2D point groups are depicted in Figure S1. and Figure S2., respectively. In the end, the matrix representation of all the group members and group generators of the 2D point groups is provided in Table S11.

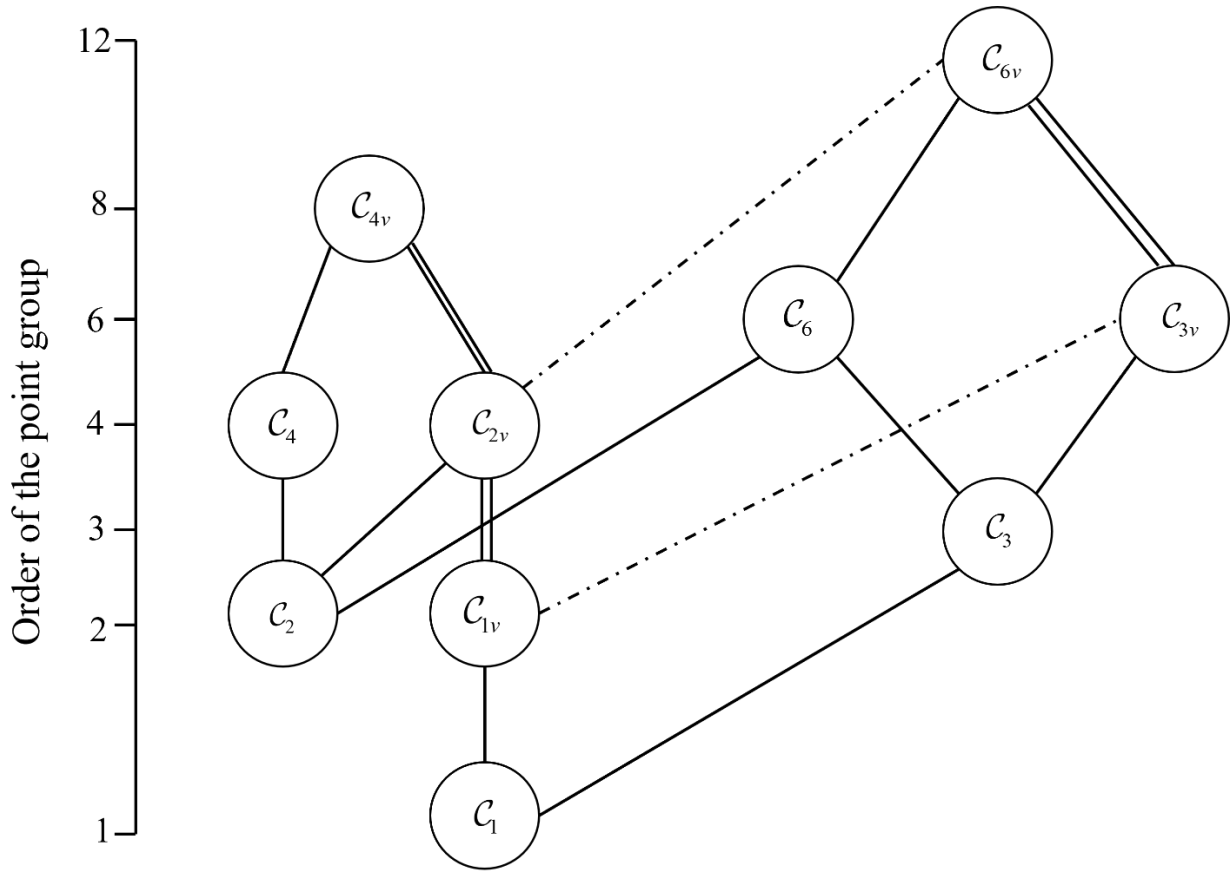

Figure S1. Symmetry hierarchy of two-dimensional crystallographic point groups.

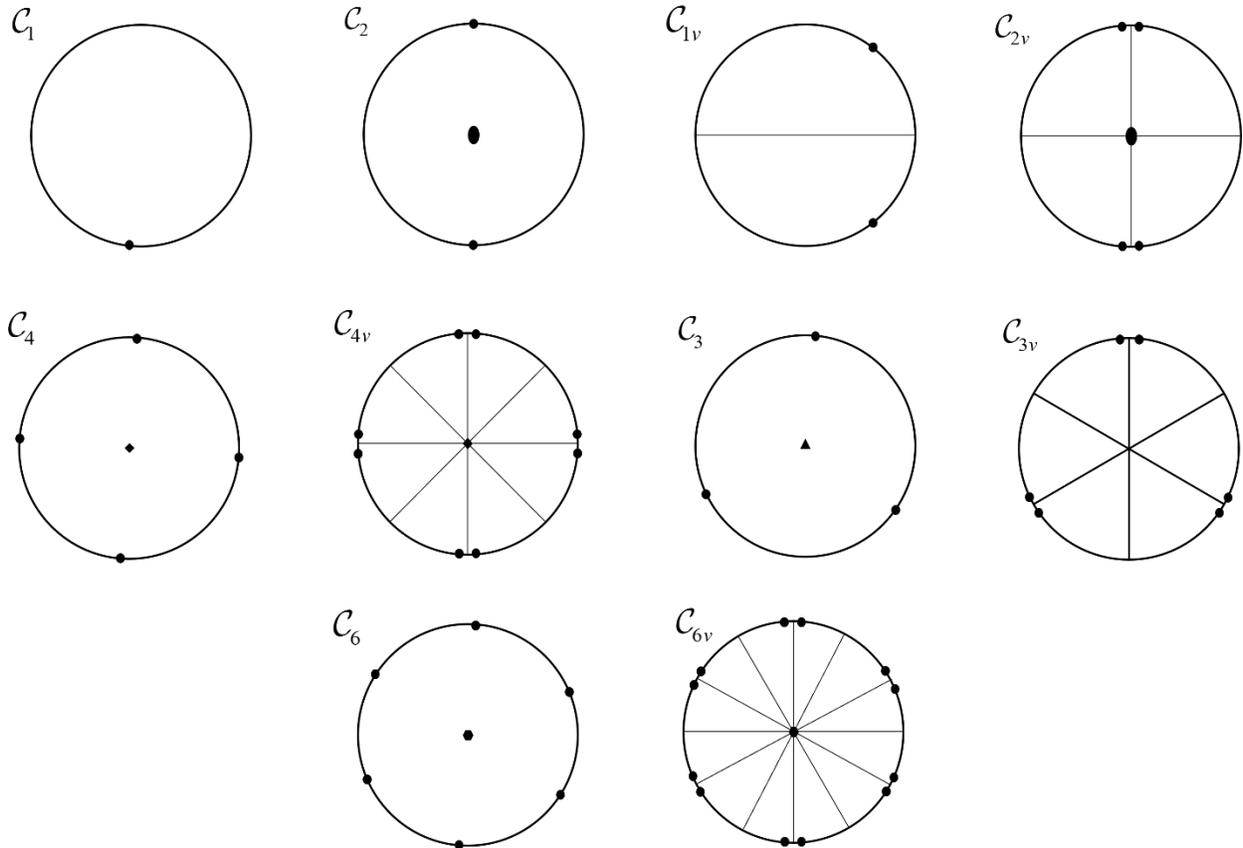

Figure S2. Graphical representation of the ten 2D crystallographic point groups.

## Multiplication table and group members of the point group $C_1$

Table S1. Multiplication table of the point group $C_1$

| Oblique($C_1$) | **E** |
|---|---|
| **E** | **E** |

Group members:
$$\mathbf{E} = \begin{pmatrix} 1 & 0 \\ 0 & 1 \end{pmatrix}$$
$\mathcal{G} = \{\mathbf{E}\}$

# Multiplication table and group members of the point group $\mathcal{C}_2$

Table S2. Multiplication table of the point group $\mathcal{C}_2$

| Oblique ($\mathcal{C}_2$) | $\mathbf{Q}^\pi$ | $\mathbf{E}$ |
|---|---|---|
| $\mathbf{Q}^\pi$ | $\mathbf{E}$ | $\mathbf{Q}^\pi$ |
| $\mathbf{E}$ | $\mathbf{Q}^\pi$ | $\mathbf{E}$ |

Group members:
$$\mathbf{Q}^\pi = \begin{pmatrix} -1 & 0 \\ 0 & -1 \end{pmatrix}, \quad \mathbf{E} = \begin{pmatrix} 1 & 0 \\ 0 & 1 \end{pmatrix}$$
$$\mathcal{G} = \{\mathbf{Q}^\pi, \mathbf{E}\}$$

# Multiplication table and group members of the point group $\mathcal{C}_{1v}$

Table S3. Multiplication table of the point group $\mathcal{C}_{1v}$

| Rectangular($\mathcal{C}_{1v}$) | $\mathbf{E}$ | $\mathbf{m}_{10}$ |
|---|---|---|
| $\mathbf{E}$ | $\mathbf{E}$ | $\mathbf{m}_{10}$ |
| $\mathbf{m}_{10}$ | $\mathbf{m}_{10}$ | $\mathbf{E}$ |

Group members:
$$\mathbf{m}_{10} = \begin{pmatrix} -1 & 0 \\ 0 & 1 \end{pmatrix}, \quad \mathbf{E} = \begin{pmatrix} 1 & 0 \\ 0 & 1 \end{pmatrix}$$
$$\mathcal{G} = \{\mathbf{E}, \mathbf{m}_{10}\}$$

# Multiplication table and group members of the point group $\mathcal{C}_{2v}$

Table S4. Multiplication table of the point group $\mathcal{C}_{2v}$

| Rectangular($\mathcal{C}_{2v}$) | $\mathbf{Q}^\pi$ | $\mathbf{E}$ | $\mathbf{m}_{10}$ | $\mathbf{m}_{01}$ |
|---|---|---|---|---|
| $\mathbf{Q}^\pi$ | $\mathbf{E}$ | $\mathbf{Q}^\pi$ | $\mathbf{m}_{01}$ | $\mathbf{m}_{10}$ |
| $\mathbf{E}$ | $\mathbf{Q}^\pi$ | $\mathbf{E}$ | $\mathbf{m}_{10}$ | $\mathbf{m}_{01}$ |
| $\mathbf{m}_{10}$ | $\mathbf{m}_{01}$ | $\mathbf{m}_{10}$ | $\mathbf{E}$ | $\mathbf{Q}^\pi$ |
| $\mathbf{m}_{01}$ | $\mathbf{m}_{10}$ | $\mathbf{m}_{01}$ | $\mathbf{Q}^\pi$ | $\mathbf{E}$ |

Group members:

$$\mathbf{m}_{10} = \begin{pmatrix} -1 & 0 \\ 0 & 1 \end{pmatrix}, \quad \mathbf{m}_{01} = \begin{pmatrix} 1 & 0 \\ 0 & -1 \end{pmatrix}, \quad \mathbf{Q}^\pi = \begin{pmatrix} -1 & 0 \\ 0 & -1 \end{pmatrix}, \quad \mathbf{E} = \begin{pmatrix} 1 & 0 \\ 0 & 1 \end{pmatrix}$$

$\mathcal{G} = \{\mathbf{E}, \mathbf{m}_{10}, \mathbf{m}_{01}, \mathbf{Q}^\pi\}$

***Remark 1***. The subscript of the mirror is the normal direction of the mirror plane.

# Multiplication table and group members of the point group $\mathcal{C}_4$

Table S5. Multiplication table of the point group $\mathcal{C}_4$

| Square($\mathcal{C}_4$) | $\mathbf{E}$ | $\mathbf{Q}^{\pi/2}$ | $\mathbf{Q}^\pi$ | $\mathbf{Q}^{3\pi/2}$ |
|---|---|---|---|---|
| $\mathbf{E}$ | $\mathbf{E}$ | $\mathbf{Q}^{\pi/2}$ | $\mathbf{Q}^\pi$ | $\mathbf{Q}^{3\pi/2}$ |
| $\mathbf{Q}^{3\pi/2}$ | $\mathbf{Q}^{3\pi/2}$ | $\mathbf{E}$ | $\mathbf{Q}^{\pi/2}$ | $\mathbf{Q}^\pi$ |
| $\mathbf{Q}^\pi$ | $\mathbf{Q}^\pi$ | $\mathbf{Q}^{3\pi/2}$ | $\mathbf{E}$ | $\mathbf{Q}^{\pi/2}$ |
| $\mathbf{Q}^{\pi/2}$ | $\mathbf{Q}^{\pi/2}$ | $\mathbf{Q}^\pi$ | $\mathbf{Q}^{3\pi/2}$ | $\mathbf{E}$ |

Group members:

$$\mathbf{Q}^{\pi/2} = \begin{pmatrix} 0 & 1 \\ -1 & 0 \end{pmatrix}, \quad \mathbf{Q}^\pi = \begin{pmatrix} -1 & 0 \\ 0 & -1 \end{pmatrix}, \quad \mathbf{Q}^{3\pi/2} = \begin{pmatrix} 0 & -1 \\ 1 & 0 \end{pmatrix}, \quad \mathbf{E} = \begin{pmatrix} 1 & 0 \\ 0 & 1 \end{pmatrix}$$

$\mathcal{G} = \{\mathbf{Q}^{\pi/2}, \mathbf{Q}^\pi, \mathbf{Q}^{3\pi/2}, \mathbf{E}\}$

# Multiplication table and group members of the point group $\mathcal{C}_{4v}$

Table S6. Multiplication table of the point group $\mathcal{C}_{4v}$

| Square($\mathcal{C}_{4v}$) | E | $Q^{\pi/2}$ | $Q^{\pi}$ | $Q^{3\pi/2}$ | $m_{10}$ | $m_{01}$ | $m_{11}$ | $m_{1-1}$ |
|---|---|---|---|---|---|---|---|---|
| E | E | $Q^{\pi/2}$ | $Q^{\pi}$ | $Q^{3\pi/2}$ | $m_{10}$ | $m_{01}$ | $m_{11}$ | $m_{1-1}$ |
| $Q^{3\pi/2}$ | $Q^{3\pi/2}$ | E | $Q^{\pi/2}$ | $Q^{\pi}$ | $m_{11}$ | $m_{1-1}$ | $m_{01}$ | $m_{10}$ |
| $Q^{\pi}$ | $Q^{\pi}$ | $Q^{3\pi/2}$ | E | $Q^{\pi/2}$ | $m_{01}$ | $m_{10}$ | $m_{1-1}$ | $m_{11}$ |
| $Q^{\pi/2}$ | $Q^{\pi/2}$ | $Q^{\pi}$ | $Q^{3\pi/2}$ | E | $m_{1-1}$ | $m_{11}$ | $m_{10}$ | $m_{01}$ |
| $m_{10}$ | $m_{10}$ | $m_{11}$ | $m_{01}$ | $m_{1-1}$ | E | $Q^{\pi}$ | $Q^{\pi/2}$ | $Q^{3\pi/2}$ |
| $m_{01}$ | $m_{01}$ | $m_{1-1}$ | $m_{10}$ | $m_{11}$ | $Q^{\pi}$ | E | $Q^{3\pi/2}$ | $Q^{\pi/2}$ |
| $m_{11}$ | $m_{11}$ | $m_{01}$ | $m_{1-1}$ | $m_{10}$ | $Q^{3\pi/2}$ | $Q^{\pi/2}$ | E | $Q^{\pi}$ |
| $m_{1-1}$ | $m_{1-1}$ | $m_{10}$ | $m_{11}$ | $m_{01}$ | $Q^{\pi/2}$ | $Q^{3\pi/2}$ | $Q^{\pi}$ | E |

Group members:

$$\mathbf{m}_{11} = \begin{pmatrix} 0 & -1 \\ -1 & 0 \end{pmatrix}, \mathbf{m}_{10} = \begin{pmatrix} -1 & 0 \\ 0 & 1 \end{pmatrix}, \mathbf{m}_{01} = \begin{pmatrix} 1 & 0 \\ 0 & -1 \end{pmatrix}, \mathbf{m}_{1-1} = \begin{pmatrix} 0 & 1 \\ 1 & 0 \end{pmatrix}, \mathbf{Q}^{\pi/2} = \begin{pmatrix} 0 & 1 \\ -1 & 0 \end{pmatrix},$$

$$\mathbf{Q}^{\pi} = \begin{pmatrix} -1 & 0 \\ 0 & -1 \end{pmatrix}, \mathbf{Q}^{3\pi/2} = \begin{pmatrix} 0 & -1 \\ 1 & 0 \end{pmatrix}, \mathbf{E} = \begin{pmatrix} 1 & 0 \\ 0 & 1 \end{pmatrix}$$

$\mathcal{G} = \{\mathbf{Q}^{\pi/2}, \mathbf{Q}^{\pi}, \mathbf{Q}^{3\pi/2}, \mathbf{E}, \mathbf{m}_{10}, \mathbf{m}_{01}, \mathbf{m}_{11}, \mathbf{m}_{1-1}\}$

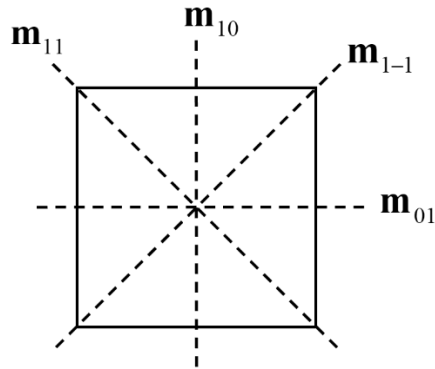

Figure S3. The representation of mirror symmetry operators for point group $\mathcal{C}_{4v}$

# Multiplication table and group members of the point group $\mathcal{C}_3$

Table S7. Multiplication table of the point group $\mathcal{C}_3$

| Trigonal($\mathcal{C}_3$) | **E** | $\mathbf{Q}^{2\pi/3}$ | $\mathbf{Q}^{4\pi/3}$ |
|---|---|---|---|
| **E** | **E** | $\mathbf{Q}^{2\pi/3}$ | $\mathbf{Q}^{4\pi/3}$ |
| $\mathbf{Q}^{4\pi/3}$ | $\mathbf{Q}^{4\pi/3}$ | **E** | $\mathbf{Q}^{2\pi/3}$ |
| $\mathbf{Q}^{2\pi/3}$ | $\mathbf{Q}^{2\pi/3}$ | $\mathbf{Q}^{4\pi/3}$ | **E** |

Group members:

$$\mathbf{Q}^{2\pi/3} = \begin{pmatrix} -1/2 & \sqrt{3}/2 \\ -\sqrt{3}/2 & -1/2 \end{pmatrix}, \quad \mathbf{Q}^{4\pi/3} = \begin{pmatrix} -1/2 & -\sqrt{3}/2 \\ \sqrt{3}/2 & -1/2 \end{pmatrix}, \quad \mathbf{E} = \begin{pmatrix} 1 & 0 \\ 0 & 1 \end{pmatrix}$$

$$\mathcal{G} = \{\mathbf{Q}^{2\pi/3}, \mathbf{Q}^{4\pi/3}, \mathbf{E}\}$$

# Multiplication table and group members of the point group $\mathcal{C}_{3v}$

Table S8. Multiplication table of the point group $\mathcal{C}_{3v}$

| Trigonal($\mathcal{C}_{3v}$) | E | $\mathbf{Q}^{2\pi/3}$ | $\mathbf{Q}^{4\pi/3}$ | $\sigma_1$ | $\sigma_2$ | $\sigma_3$ |
|---|---|---|---|---|---|---|
| E | E | $\mathbf{Q}^{2\pi/3}$ | $\mathbf{Q}^{4\pi/3}$ | $\sigma_1$ | $\sigma_2$ | $\sigma_3$ |
| $\mathbf{Q}^{4\pi/3}$ | $\mathbf{Q}^{4\pi/3}$ | E | $\mathbf{Q}^{2\pi/3}$ | $\sigma_2$ | $\sigma_3$ | $\sigma_1$ |
| $\mathbf{Q}^{2\pi/3}$ | $\mathbf{Q}^{2\pi/3}$ | $\mathbf{Q}^{4\pi/3}$ | E | $\sigma_3$ | $\sigma_1$ | $\sigma_2$ |
| $\sigma_1$ | $\sigma_1$ | $\sigma_2$ | $\sigma_3$ | E | $\mathbf{Q}^{2\pi/3}$ | $\mathbf{Q}^{4\pi/3}$ |
| $\sigma_2$ | $\sigma_2$ | $\sigma_3$ | $\sigma_1$ | $\mathbf{Q}^{4\pi/3}$ | E | $\mathbf{Q}^{2\pi/3}$ |
| $\sigma_3$ | $\sigma_3$ | $\sigma_1$ | $\sigma_2$ | $\mathbf{Q}^{2\pi/3}$ | $\mathbf{Q}^{4\pi/3}$ | E |

Group members:

$$\mathbf{Q}^{2\pi/3} = \begin{pmatrix} -1/2 & \sqrt{3}/2 \\ -\sqrt{3}/2 & -1/2 \end{pmatrix},\ \mathbf{Q}^{4\pi/3} = \begin{pmatrix} -1/2 & -\sqrt{3}/2 \\ \sqrt{3}/2 & -1/2 \end{pmatrix},\ \mathbf{E} = \begin{pmatrix} 1 & 0 \\ 0 & 1 \end{pmatrix},\ \boldsymbol{\sigma}_1 = \begin{pmatrix} -1 & 0 \\ 0 & 1 \end{pmatrix},$$

$$\boldsymbol{\sigma}_2 = \begin{pmatrix} 1/2 & -\sqrt{3}/2 \\ -\sqrt{3}/2 & -1/2 \end{pmatrix},\ \boldsymbol{\sigma}_3 = \begin{pmatrix} 1/2 & \sqrt{3}/2 \\ \sqrt{3}/2 & -1/2 \end{pmatrix}$$

$$\mathcal{G} = \{\mathbf{Q}^{2\pi/3}, \mathbf{Q}^{4\pi/3}, \mathbf{E}, \boldsymbol{\sigma}_1, \boldsymbol{\sigma}_2, \boldsymbol{\sigma}_3\}$$

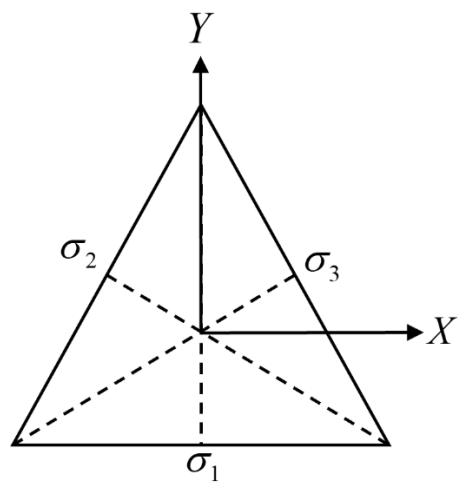

Figure S4. The representation of mirror symmetry operators for point group $\mathcal{C}_{3v}$

# Multiplication table and group members of the point group $\mathcal{C}_6$

Table S9. Multiplication table of the point group $\mathcal{C}_6$

| Hexagonal($\mathcal{C}_6$) | **E** | $\mathbf{Q}^{\pi/3}$ | $\mathbf{Q}^{2\pi/3}$ | $\mathbf{Q}^{\pi}$ | $\mathbf{Q}^{4\pi/3}$ | $\mathbf{Q}^{5\pi/3}$ |
|---|---|---|---|---|---|---|
| **E** | **E** | $\mathbf{Q}^{\pi/3}$ | $\mathbf{Q}^{2\pi/3}$ | $\mathbf{Q}^{\pi}$ | $\mathbf{Q}^{4\pi/3}$ | $\mathbf{Q}^{5\pi/3}$ |
| $\mathbf{Q}^{5\pi/3}$ | $\mathbf{Q}^{5\pi/3}$ | **E** | $\mathbf{Q}^{\pi/3}$ | $\mathbf{Q}^{2\pi/3}$ | $\mathbf{Q}^{\pi}$ | $\mathbf{Q}^{4\pi/3}$ |
| $\mathbf{Q}^{4\pi/3}$ | $\mathbf{Q}^{4\pi/3}$ | $\mathbf{Q}^{5\pi/3}$ | **E** | $\mathbf{Q}^{\pi/3}$ | $\mathbf{Q}^{2\pi/3}$ | $\mathbf{Q}^{\pi}$ |
| $\mathbf{Q}^{\pi}$ | $\mathbf{Q}^{\pi}$ | $\mathbf{Q}^{4\pi/3}$ | $\mathbf{Q}^{5\pi/3}$ | **E** | $\mathbf{Q}^{\pi/3}$ | $\mathbf{Q}^{2\pi/3}$ |
| $\mathbf{Q}^{2\pi/3}$ | $\mathbf{Q}^{2\pi/3}$ | $\mathbf{Q}^{\pi}$ | $\mathbf{Q}^{4\pi/3}$ | $\mathbf{Q}^{5\pi/3}$ | **E** | $\mathbf{Q}^{\pi/3}$ |
| $\mathbf{Q}^{\pi/3}$ | $\mathbf{Q}^{\pi/3}$ | $\mathbf{Q}^{2\pi/3}$ | $\mathbf{Q}^{\pi}$ | $\mathbf{Q}^{4\pi/3}$ | $\mathbf{Q}^{5\pi/3}$ | **E** |

Group members:

$$\mathbf{Q}^{\pi/3} = \begin{pmatrix} 1/2 & \sqrt{3}/2 \\ -\sqrt{3}/2 & 1/2 \end{pmatrix}, \mathbf{Q}^{2\pi/3} = \begin{pmatrix} -1/2 & \sqrt{3}/2 \\ -\sqrt{3}/2 & -1/2 \end{pmatrix}, \mathbf{Q}^{\pi} = \begin{pmatrix} -1 & 0 \\ 0 & -1 \end{pmatrix},$$

$$\mathbf{Q}^{4\pi/3} = \begin{pmatrix} -1/2 & -\sqrt{3}/2 \\ \sqrt{3}/2 & -1/2 \end{pmatrix}, \mathbf{Q}^{5\pi/3} = \begin{pmatrix} 1/2 & -\sqrt{3}/2 \\ \sqrt{3}/2 & 1/2 \end{pmatrix}, \mathbf{E} = \begin{pmatrix} 1 & 0 \\ 0 & 1 \end{pmatrix}$$

$\mathcal{G} = \{\mathbf{Q}^{\pi/3}, \mathbf{Q}^{2\pi/3}, \mathbf{Q}^{\pi}, \mathbf{Q}^{4\pi/3}, \mathbf{Q}^{5\pi/3}, \mathbf{E}\}$

# Multiplication table and group members of the point group $\mathcal{C}_{6v}$

Table S10. Multiplication table of the point group $\mathcal{C}_{6v}$

| Hexagonal ($\mathcal{C}_{6v}$) | E | $Q^{\pi/3}$ | $Q^{2\pi/3}$ | $Q^{\pi}$ | $Q^{4\pi/3}$ | $Q^{5\pi/3}$ | $\sigma_{d_1}$ | $\sigma_{d_2}$ | $\sigma_{d_3}$ | $\sigma_{v_1}$ | $\sigma_{v_2}$ | $\sigma_{v_3}$ |
|---|---|---|---|---|---|---|---|---|---|---|---|---|
| E | E | $Q^{\pi/3}$ | $Q^{2\pi/3}$ | $Q^{\pi}$ | $Q^{4\pi/3}$ | $Q^{5\pi/3}$ | $\sigma_{d_1}$ | $\sigma_{d_2}$ | $\sigma_{d_3}$ | $\sigma_{v_1}$ | $\sigma_{v_2}$ | $\sigma_{v_3}$ |
| $Q^{5\pi/3}$ | $Q^{5\pi/3}$ | E | $Q^{\pi/3}$ | $Q^{2\pi/3}$ | $Q^{\pi}$ | $Q^{4\pi/3}$ | $\sigma_{v_2}$ | $\sigma_{v_3}$ | $\sigma_{v_1}$ | $\sigma_{d_2}$ | $\sigma_{d_3}$ | $\sigma_{d_1}$ |
| $Q^{4\pi/3}$ | $Q^{4\pi/3}$ | $Q^{5\pi/3}$ | E | $Q^{\pi/3}$ | $Q^{2\pi/3}$ | $Q^{\pi}$ | $\sigma_{d_3}$ | $\sigma_{d_1}$ | $\sigma_{d_2}$ | $\sigma_{v_3}$ | $\sigma_{v_1}$ | $\sigma_{v_2}$ |
| $Q^{\pi}$ | $Q^{\pi}$ | $Q^{4\pi/3}$ | $Q^{5\pi/3}$ | E | $Q^{\pi/3}$ | $Q^{2\pi/3}$ | $\sigma_{v_1}$ | $\sigma_{v_2}$ | $\sigma_{v_3}$ | $\sigma_{d_1}$ | $\sigma_{d_2}$ | $\sigma_{d_3}$ |
| $Q^{2\pi/3}$ | $Q^{2\pi/3}$ | $Q^{\pi}$ | $Q^{4\pi/3}$ | $Q^{5\pi/3}$ | E | $Q^{\pi/3}$ | $\sigma_{d_2}$ | $\sigma_{d_3}$ | $\sigma_{d_1}$ | $\sigma_{v_2}$ | $\sigma_{v_3}$ | $\sigma_{v_1}$ |
| $Q^{\pi/3}$ | $Q^{\pi/3}$ | $Q^{2\pi/3}$ | $Q^{\pi}$ | $Q^{4\pi/3}$ | $Q^{5\pi/3}$ | E | $\sigma_{v_3}$ | $\sigma_{v_1}$ | $\sigma_{v_2}$ | $\sigma_{d_3}$ | $\sigma_{d_1}$ | $\sigma_{d_2}$ |
| $\sigma_{d_1}$ | $\sigma_{d_1}$ | $\sigma_{v_2}$ | $\sigma_{d_3}$ | $\sigma_{v_1}$ | $\sigma_{d_2}$ | $\sigma_{v_3}$ | E | $Q^{4\pi/3}$ | $Q^{2\pi/3}$ | $Q^{\pi}$ | $Q^{\pi/3}$ | $Q^{5\pi/3}$ |
| $\sigma_{d_2}$ | $\sigma_{d_2}$ | $\sigma_{v_3}$ | $\sigma_{d_1}$ | $\sigma_{v_2}$ | $\sigma_{d_3}$ | $\sigma_{v_1}$ | $Q^{2\pi/3}$ | E | $Q^{4\pi/3}$ | $Q^{5\pi/3}$ | $Q^{\pi}$ | $Q^{\pi/3}$ |
| $\sigma_{d_3}$ | $\sigma_{d_3}$ | $\sigma_{v_1}$ | $\sigma_{d_2}$ | $\sigma_{v_3}$ | $\sigma_{d_1}$ | $\sigma_{v_2}$ | $Q^{4\pi/3}$ | $Q^{2\pi/3}$ | E | $Q^{\pi/3}$ | $Q^{5\pi/3}$ | $Q^{\pi}$ |
| $\sigma_{v_1}$ | $\sigma_{v_1}$ | $\sigma_{d_2}$ | $\sigma_{v_3}$ | $\sigma_{d_1}$ | $\sigma_{v_2}$ | $\sigma_{d_3}$ | $Q^{\pi}$ | $Q^{\pi/3}$ | $Q^{5\pi/3}$ | E | $Q^{4\pi/3}$ | $Q^{2\pi/3}$ |
| $\sigma_{v_2}$ | $\sigma_{v_2}$ | $\sigma_{d_3}$ | $\sigma_{v_1}$ | $\sigma_{d_2}$ | $\sigma_{v_3}$ | $\sigma_{d_1}$ | $Q^{5\pi/3}$ | $Q^{\pi}$ | $Q^{\pi/3}$ | $Q^{2\pi/3}$ | E | $Q^{4\pi/3}$ |
| $\sigma_{v_3}$ | $\sigma_{v_3}$ | $\sigma_{d_1}$ | $\sigma_{v_2}$ | $\sigma_{d_3}$ | $\sigma_{v_1}$ | $\sigma_{d_2}$ | $Q^{\pi/3}$ | $Q^{5\pi/3}$ | $Q^{\pi}$ | $Q^{4\pi/3}$ | $Q^{2\pi/3}$ | E |

Group members:

$$\mathbf{Q}^{\pi/3} = \begin{pmatrix} 1/2 & \sqrt{3}/2 \\ -\sqrt{3}/2 & 1/2 \end{pmatrix}, \quad \mathbf{Q}^{2\pi/3} = \begin{pmatrix} -1/2 & \sqrt{3}/2 \\ -\sqrt{3}/2 & -1/2 \end{pmatrix}, \quad \mathbf{Q}^{\pi} = \begin{pmatrix} -1 & 0 \\ 0 & -1 \end{pmatrix},$$

$$\mathbf{Q}^{4\pi/3} = \begin{pmatrix} -1/2 & -\sqrt{3}/2 \\ \sqrt{3}/2 & -1/2 \end{pmatrix}, \quad \mathbf{Q}^{5\pi/3} = \begin{pmatrix} 1/2 & -\sqrt{3}/2 \\ \sqrt{3}/2 & 1/2 \end{pmatrix}, \quad \mathbf{E} = \begin{pmatrix} 1 & 0 \\ 0 & 1 \end{pmatrix}, \quad \boldsymbol{\sigma}_{d_1} = \begin{pmatrix} -1 & 0 \\ 0 & 1 \end{pmatrix},$$

$$\boldsymbol{\sigma}_{d_2} = \begin{pmatrix} 1/2 & \sqrt{3}/2 \\ \sqrt{3}/2 & -1/2 \end{pmatrix}, \quad \boldsymbol{\sigma}_{d_3} = \begin{pmatrix} 1/2 & -\sqrt{3}/2 \\ -\sqrt{3}/2 & -1/2 \end{pmatrix}, \quad \boldsymbol{\sigma}_{v_1} = \begin{pmatrix} 1 & 0 \\ 0 & -1 \end{pmatrix}, \quad \boldsymbol{\sigma}_{v_2} = \begin{pmatrix} -1/2 & -\sqrt{3}/2 \\ -\sqrt{3}/2 & 1/2 \end{pmatrix},$$

$$\boldsymbol{\sigma}_{v_3} = \begin{pmatrix} -1/2 & \sqrt{3}/2 \\ \sqrt{3}/2 & 1/2 \end{pmatrix}$$

$$\mathcal{G} = \{\mathbf{Q}^{\pi/3}, \mathbf{Q}^{2\pi/3}, \mathbf{Q}^{\pi}, \mathbf{Q}^{4\pi/3}, \mathbf{Q}^{5\pi/3}, \mathbf{E}, \boldsymbol{\sigma}_{d_1}, \boldsymbol{\sigma}_{d_2}, \boldsymbol{\sigma}_{d_3}, \boldsymbol{\sigma}_{v_1}, \boldsymbol{\sigma}_{v_2}, \boldsymbol{\sigma}_{v_3}\}$$

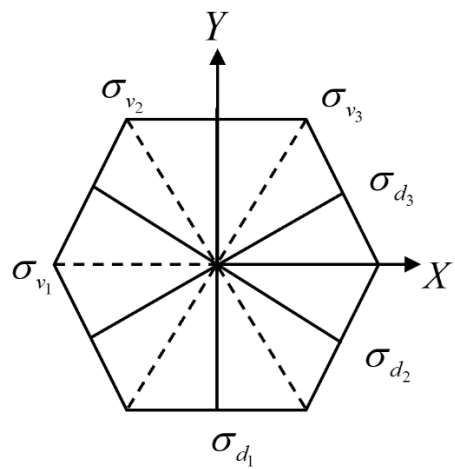

Figure S5. The representation of mirror symmetry operators for the point group $\mathcal{C}_{6v}$

Table S11. Members and generators for the two-dimensional crystallographic and compact point groups

| Lattice structure | Schoenflies Symbol (Hermann–Mauguin) | Group members | Group generators |
|---|---|---|---|
| Oblique | $C_1$ (1) | $\mathbf{E} = \begin{pmatrix} 1 & 0 \\ 0 & 1 \end{pmatrix}$ | $\mathbf{E} = \begin{pmatrix} 1 & 0 \\ 0 & 1 \end{pmatrix}$ |
| | $C_2$ (2) | $\mathbf{Q}^\pi = \begin{pmatrix} -1 & 0 \\ 0 & -1 \end{pmatrix},\ \mathbf{E} = \begin{pmatrix} 1 & 0 \\ 0 & 1 \end{pmatrix}$ | $\mathbf{Q}^\pi = \begin{pmatrix} -1 & 0 \\ 0 & -1 \end{pmatrix}$ |
| Rectangular | $C_{1v}$ (m) | $\mathbf{m}_{01} = \begin{pmatrix} 1 & 0 \\ 0 & -1 \end{pmatrix},\ \mathbf{E} = \begin{pmatrix} 1 & 0 \\ 0 & 1 \end{pmatrix}$ | $\mathbf{m}_{01} = \begin{pmatrix} 1 & 0 \\ 0 & -1 \end{pmatrix}$ |
| | $C_{2v}$ (2mm) | $\mathbf{m}_{10} = \begin{pmatrix} -1 & 0 \\ 0 & 1 \end{pmatrix},\ \mathbf{m}_{01} = \begin{pmatrix} 1 & 0 \\ 0 & -1 \end{pmatrix},$ $\mathbf{Q}^\pi = \begin{pmatrix} -1 & 0 \\ 0 & -1 \end{pmatrix},\ \mathbf{E} = \begin{pmatrix} 1 & 0 \\ 0 & 1 \end{pmatrix}$ | $\mathbf{Q}^\pi = \begin{pmatrix} -1 & 0 \\ 0 & -1 \end{pmatrix},\ \mathbf{m}_{01} = \begin{pmatrix} 1 & 0 \\ 0 & -1 \end{pmatrix}$ |
| Square | $C_4$ (4) | $\mathbf{Q}^{\pi/2} = \begin{pmatrix} 0 & 1 \\ -1 & 0 \end{pmatrix},\ \mathbf{Q}^\pi = \begin{pmatrix} -1 & 0 \\ 0 & -1 \end{pmatrix},$ $\mathbf{Q}^{3\pi/2} = \begin{pmatrix} 0 & -1 \\ 1 & 0 \end{pmatrix},\ \mathbf{E} = \begin{pmatrix} 1 & 0 \\ 0 & 1 \end{pmatrix}$ | $\mathbf{Q}^{\pi/2} = \begin{pmatrix} 0 & 1 \\ -1 & 0 \end{pmatrix}$ |

| | $\mathcal{C}_{4v}$ (4mm) | $\mathbf{m}_{11} = \begin{pmatrix} 0 & -1 \\ -1 & 0 \end{pmatrix}$, $\mathbf{m}_{10} = \begin{pmatrix} -1 & 0 \\ 0 & 1 \end{pmatrix}$, $\mathbf{m}_{01} = \begin{pmatrix} 1 & 0 \\ 0 & -1 \end{pmatrix}$, $\mathbf{m}_{1-1} = \begin{pmatrix} 0 & 1 \\ 1 & 0 \end{pmatrix}$, $\mathbf{Q}^{\pi/2} = \begin{pmatrix} 0 & 1 \\ -1 & 0 \end{pmatrix}$, $\mathbf{Q}^{\pi} = \begin{pmatrix} -1 & 0 \\ 0 & -1 \end{pmatrix}$, $\mathbf{Q}^{3\pi/2} = \begin{pmatrix} 0 & -1 \\ 1 & 0 \end{pmatrix}$, $\mathbf{E} = \begin{pmatrix} 1 & 0 \\ 0 & 1 \end{pmatrix}$ | $\mathbf{m}_{10} = \begin{pmatrix} -1 & 0 \\ 0 & 1 \end{pmatrix}$, $\mathbf{Q}^{\pi/2} = \begin{pmatrix} 0 & 1 \\ -1 & 0 \end{pmatrix}$ |
|---|---|---|---|
| Trigonal | $\mathcal{C}_3$ (3) | $\mathbf{Q}^{2\pi/3} = \begin{pmatrix} -1/2 & \sqrt{3}/2 \\ -\sqrt{3}/2 & -1/2 \end{pmatrix}$, $\mathbf{Q}^{4\pi/3} = \begin{pmatrix} -1/2 & -\sqrt{3}/2 \\ \sqrt{3}/2 & -1/2 \end{pmatrix}$, $\mathbf{E} = \begin{pmatrix} 1 & 0 \\ 0 & 1 \end{pmatrix}$ | $\mathbf{Q}^{2\pi/3} = \begin{pmatrix} -1/2 & \sqrt{3}/2 \\ -\sqrt{3}/2 & -1/2 \end{pmatrix}$ |
| | $\mathcal{C}_{3v}$ (3m) | $\mathbf{Q}^{2\pi/3} = \begin{pmatrix} -1/2 & \sqrt{3}/2 \\ -\sqrt{3}/2 & -1/2 \end{pmatrix}$, $\mathbf{Q}^{4\pi/3} = \begin{pmatrix} -1/2 & -\sqrt{3}/2 \\ \sqrt{3}/2 & -1/2 \end{pmatrix}$, $\mathbf{E} = \begin{pmatrix} 1 & 0 \\ 0 & 1 \end{pmatrix}$, $\boldsymbol{\sigma}_1 = \begin{pmatrix} -1 & 0 \\ 0 & 1 \end{pmatrix}$, $\boldsymbol{\sigma}_2 = \begin{pmatrix} 1/2 & -\sqrt{3}/2 \\ -\sqrt{3}/2 & -1/2 \end{pmatrix}$, $\boldsymbol{\sigma}_3 = \begin{pmatrix} 1/2 & \sqrt{3}/2 \\ \sqrt{3}/2 & -1/2 \end{pmatrix}$ | $\mathbf{Q}^{2\pi/3} = \begin{pmatrix} -1/2 & \sqrt{3}/2 \\ -\sqrt{3}/2 & -1/2 \end{pmatrix}$, $\boldsymbol{\sigma}_1 = \begin{pmatrix} -1 & 0 \\ 0 & 1 \end{pmatrix}$ |

| | | | |
|---|---|---|---|
| Hexagonal | $\mathcal{C}_6$ (6) | $\mathbf{Q}^{\pi/3} = \begin{pmatrix} 1/2 & \sqrt{3}/2 \\ -\sqrt{3}/2 & 1/2 \end{pmatrix},$ $\mathbf{Q}^{2\pi/3} = \begin{pmatrix} -1/2 & \sqrt{3}/2 \\ -\sqrt{3}/2 & -1/2 \end{pmatrix},$ $\mathbf{Q}^{\pi} = \begin{pmatrix} -1 & 0 \\ 0 & -1 \end{pmatrix},$ $\mathbf{Q}^{4\pi/3} = \begin{pmatrix} -1/2 & -\sqrt{3}/2 \\ \sqrt{3}/2 & -1/2 \end{pmatrix},$ $\mathbf{Q}^{5\pi/3} = \begin{pmatrix} 1/2 & -\sqrt{3}/2 \\ \sqrt{3}/2 & 1/2 \end{pmatrix},$ $\mathbf{E} = \begin{pmatrix} 1 & 0 \\ 0 & 1 \end{pmatrix}$ | $\mathbf{Q}^{\pi/3} = \begin{pmatrix} 1/2 & \sqrt{3}/2 \\ -\sqrt{3}/2 & 1/2 \end{pmatrix}$ |
| | $\mathcal{C}_{6v}$ (6mm) | $\mathbf{Q}^{\pi/3} = \begin{pmatrix} 1/2 & \sqrt{3}/2 \\ -\sqrt{3}/2 & 1/2 \end{pmatrix},$ $\mathbf{Q}^{2\pi/3} = \begin{pmatrix} -1/2 & \sqrt{3}/2 \\ -\sqrt{3}/2 & -1/2 \end{pmatrix},$ $\mathbf{Q}^{\pi} = \begin{pmatrix} -1 & 0 \\ 0 & -1 \end{pmatrix},$ $\mathbf{Q}^{4\pi/3} = \begin{pmatrix} -1/2 & -\sqrt{3}/2 \\ \sqrt{3}/2 & -1/2 \end{pmatrix},$ $\mathbf{Q}^{5\pi/3} = \begin{pmatrix} 1/2 & -\sqrt{3}/2 \\ \sqrt{3}/2 & 1/2 \end{pmatrix},$ $\mathbf{E} = \begin{pmatrix} 1 & 0 \\ 0 & 1 \end{pmatrix},$ $\boldsymbol{\sigma}_{d_1} = \begin{pmatrix} -1 & 0 \\ 0 & 1 \end{pmatrix},$ $\boldsymbol{\sigma}_{d_2} = \begin{pmatrix} 1/2 & \sqrt{3}/2 \\ \sqrt{3}/2 & -1/2 \end{pmatrix},$ $\boldsymbol{\sigma}_{d_3} = \begin{pmatrix} 1/2 & -\sqrt{3}/2 \\ -\sqrt{3}/2 & -1/2 \end{pmatrix},$ $\boldsymbol{\sigma}_{v_1} = \begin{pmatrix} 1 & 0 \\ 0 & -1 \end{pmatrix},$ $\boldsymbol{\sigma}_{v_2} = \begin{pmatrix} -1/2 & -\sqrt{3}/2 \\ -\sqrt{3}/2 & 1/2 \end{pmatrix},$ $\boldsymbol{\sigma}_{v_3} = \begin{pmatrix} -1/2 & \sqrt{3}/2 \\ \sqrt{3}/2 & 1/2 \end{pmatrix}$ | $\mathbf{Q}^{\pi/3} = \begin{pmatrix} 1/2 & \sqrt{3}/2 \\ -\sqrt{3}/2 & 1/2 \end{pmatrix},$ $\boldsymbol{\sigma}_{d_1} = \begin{pmatrix} -1 & 0 \\ 0 & 1 \end{pmatrix}$ |

| Circular | $\mathcal{C}_\infty$ | Infinite number of members | $\mathbf{Q}^\theta = \begin{pmatrix} \cos\theta & \sin\theta \\ -\sin\theta & \cos\theta \end{pmatrix}$ |
| --- | --- | --- | --- |
| | $\mathcal{C}_{\infty v}$ | Infinite number of members | $\mathbf{Q}^\theta = \begin{pmatrix} \cos\theta & \sin\theta \\ -\sin\theta & \cos\theta \end{pmatrix}, \boldsymbol{\sigma}_1 = \begin{pmatrix} -1 & 0 \\ 0 & 1 \end{pmatrix}$ |



\end{document}